\documentclass[]{article}

\usepackage{amsfonts,amsmath,amssymb}
\usepackage{amsthm}
\usepackage[cp866]{inputenc}
\usepackage[english]{babel}

\allowdisplaybreaks \theoremstyle{definition}
\newtheorem{Definition}{Definition}

\theoremstyle{plain}
\newtheorem{Theorem}{Theorem}
\newtheorem{Corollary}{Corollary}
\newtheorem{Proposition}{Proposition}

\newtheorem{Lemma}{Lemma}

\newcommand{\U}{\mathcal{U}}
\newcommand{\con}{\natural}

\newcommand{\Pomega}{{P\hspace*{-1pt}\omega}}
\newcommand{\bfSig}{\mathbf{\Sigma}}
\newcommand{\bfPi}{\mathbf{\Pi}}
\newcommand{\bfDelta}{\mathbf{\Delta}}
\newcommand{\bfGamma}{\mathbf{\Gamma}}

\newcommand{\calN}{\mathcal{N}}

\newcommand{\IR}{\mathbb{R}}
\title{A $Q$-Wadge Hierarchy in Quasi-Polish Spaces}
\author{Victor Selivanov
\\A.P. Ershov
Institute of
Informatics Systems SB RAS\\ 
{\tt vseliv@iis.nsk.su}
}

\begin{document}
\large
\date{}
 \maketitle

\begin{abstract}
The Wadge hierarchy was originally defined and studied only in the Baire space (and some other zero-dimensional spaces). We extend it here to arbitrary topological spaces by providing a set-theoretic definition of all its levels. We show that our extension behaves well in  second countable spaces and especially in quasi-Polish spaces. In particular, all levels are preserved by continuous open surjections between second countable spaces which implies e.g. several Hausdorff-Kuratowski-type theorems in quasi-Polish spaces. In fact, many results hold not only for the Wadge hierarchy of sets but also for its extension to Borel functions from a space to a countable better quasiorder $Q$. 

 {\bf Key words}: Borel hierarchy, Wadge hierarchy, fine hierarchy, iterated labeled tree, $h$-quasiorder, better quasiorder, $Q$-partition.

%{\bf MSC}: Primary   03D15, 03D78, 58J45; Secondary 65M06, 65M25.
\end{abstract}

\section{Introduction}\label{in}

The classical Borel, Luzin, and Hausdorff hierarchies in Polish spaces, which are defined  using set operations, play an important role in descriptive set theory (DST). Recently, these hierarchies were extended and shown to have similar nice properties also in quasi-Polish spaces \cite{br} which include many non-Hausdorff spaces of interest for several branches of mathematics and theoretical computer science.

The Wadge hierarchy, introduced in \cite{wad72,wad84}, is  non-classical in the sense that it is based on a notion of reducibility that was not recognized in the classical DST, and on using ingenious versions of Gale-Stewart games rather than on  set operations.
For subsets $A,B$ of the Baire space $\mathcal{N}=\omega^\omega$, $A$ is {\em Wadge reducible} to $B$ ($A\leq_WB$), if $A=f^{-1}(B)$ for some continuous function $f$ on $\mathcal{N}$. The quotient-poset of the preorder $(P(\mathcal{N});\leq_W)$ under the induced equivalence relation $\equiv_W$ on the power-set of $\mathcal{N}$ is called {\em the structure of Wadge degrees} in $\mathcal{N}$. W. Wadge \cite{wad84}  characterised the structure of Wadge degrees of Borel sets (i.e., the quotient-poset of $({\mathbf B}(\mathcal{N});\leq_W)$) up to isomorphism. In particular, this quotient-poset is semi-well-ordered, hence it is well-founded and has no 3 pairwise incomparable elements. For more information on Wadge degrees see \cite{vw76,kls12}. 

This gives rise to the Wadge hierarchy $\{\bfSig_\alpha(\calN)\}_{\alpha<\upsilon}$ (for a rather large ordinal $\upsilon$) in $\calN$ which is a great refinement of the Borel hierarchy (for more information see the next section where we also give precise definitions of other notions mentioned in this introduction).
The Wadge hierarchy was originally defined only for the Baire space. Using the methods of \cite{wad84} it is easy to check that the structure $({\mathbf B}(X);\leq_W)$ of Wadge degrees of Borel sets in any zero-dimensional Polish space $X$ remains semi-well-ordered and the Wadge hierarchy in such spaces looks rather similar to that in the Baire space.

The Wadge hierarchy of sets was an important development in classical DST not only as a unifying concept (it subsumes all hierarchies known before) but also as a useful tool to investigate second countable zero-dimension spaces. We illustrate this with two examples. In \cite{ems87}  a complete classification (up to homeomorphism) of homogeneous zero-dimensional absolute Borel sets was achieved, completing a series of earlier results in this direction. In \cite{ems87} it was shown that any Borel subspace of the Baire space with more than one point has a non-trivial auto-homeomorphism.

In this paper we attempt to find the ``correct'' extension of the Wadge hierarchy from Polish zero-dimensional spaces to arbitrary second countable spaces. There are at least three approaches to this problem.

The first approach is to show that Wadge reducibility in such spaces behaves similarly to its behaviour in the Baire space, i.e. it is a semi-well-order. Unfortunately, this is not the case: for many natural quasi-Polish spaces $X$ the structure $({\mathbf B}(X);\leq_W)$ is not well-founded and has antichains with more than 2 elements (see e.g. \cite{he96,sc10,s05,ik10,bg15,du19}). Thus, this approach does not lead to a reasonable extension of the Wadge hierarchy to quasi-Polish spaces.

The second  approach  was independently suggested in \cite{pe15,s16}. The approach is based on the characterization of quasi-Polish spaces as the second countable $T_0$-spaces $X$ such that there is a total admissible representation $\xi$ from $\calN$ onto $X$ \cite{br}. Namely, one can  {\em define} the Wadge hierarchy $\{\bfSig_\alpha(X)\}_{\alpha<\upsilon}$ in $X$ by $\bfSig_\alpha(X)=\{A\subseteq X\mid \xi^{-1}(A)\in\bfSig_\alpha(\calN)\}$. One easily checks that the definition of $\bfSig_\alpha(X)$ does not depend on the choice of $\xi$, $\bigcup_{\alpha<\nu}\bfSig_\alpha(X)=\bf{B}(X)$,   ${\bfSig}_\alpha(X)\subseteq{\bfDelta}_\beta(X)$ for all $\alpha<\beta<\upsilon$, and any $\bfSig_\alpha(X)$ is downward closed under the Wadge reducibility in $X$. This definition is short and elegant but it gives no real understanding of how the levels $\bfSig_\alpha(X)$ look like, in particular their set-theoretic descriptions are completely unclear. 

The third  approach (traditional in classical DST) proposed in \cite{s16} is to apply a refinement process according to which one starts with the Borel hierarchy and subsequently defines suitable ``natural'' refinements of the hierarchies already available. At the first step of this process we obtain the  Hausdorff hierarchies over each level of the Borel hierarchy thoroughly investigated in \cite{br}. Further refinements may be done using more sophisticated set operations which extend and modify some operations introduced in \cite{wad84} for the Baire space. In this way we described in \cite{s16,s17} an increasing sequence of pointclasses $\{\bfSig_\alpha(X)\}_{\alpha<\lambda}$, $\lambda=sup\{\omega_1,\omega_1^{\omega_1},\omega_1^{\omega_1^{\omega_1}},\ldots\}$ which exhaust the sets of finite Borel rank, and we conjectured their coincidence with the corresponding classes from the second approach and proved the conjecture for some particular cases. Thus, we proposed a way to achieving a reasonable set-theoretic definition of the  Wadge hierarchy in $X$.

In the present paper we propose such a definition for the whole Wadge hierarchy. The definition  is an infinitary version of  the so called fine hierarchy  introduced and studied in a series of my publications (see e.g. \cite{s08} for a survey). In fact, this paper develops a ``classical''  infinitary version of the effective finitary version of the Wadge hierarchy in effective spaces and computable quasi-Polish spaces recently developed in \cite{s19}. Arguably, our infinitary fine hierarchy (IFH), and hence also the Wadge hierarchy, is a kind of ``iterated difference hierarchy'' over levels of the Borel hierarchy; it only remains to make precise how to ``iterate'' the difference hierarchies.

 Along with describing (hopefully) the right version of the Wadge hierarchy (by identifying it with the IFH) in arbitrary spaces we show that this version behaves well in  second countable spaces and especially in quasi-Polish spaces. E.g., it provides the description  of all levels $\bfSig_\alpha(X)$ in quasi-Polish spaces. Also, all levels of the IFH are preserved by continuous open surjections between second countable spaces which gives a broad extension of results by Saint Raymond and de Brecht for the Borel and Hausdorff hierarchies  \cite{sr07,br}. We also show that several Hausdorff-Kuratowski-type theorems are inherited by the continuous open images. As a corollary we obtain  such theorems in arbitrary quasi-Polish spaces. %Our approach also yields extensions and new proofs of some known results, in particular the Steel separation theorem for the non-self-dual levels of the Wadge hierarchy \cite{ste80} becomes a simple corollary of the reduction property which again implies some generalisations.  
 
 The notions and results of this paper apply not only to the Wadge hierarchy of sets  discussed so far but also to a more general hierarchy of  functions $A:X\to Q$ from a space $X$ to an arbitrary quasiorder $Q$. We identify such functions with {\em $Q$-partitions of $X$} of the form $\{A^{-1}(q)\}_{q\in Q}$ in order to stress  their close relation  to $k$-partitions (obtained when $Q=\bar{k}=\{0,\dots,k-1\}$ is an antichain with $k$-elements) studied by several authors. 

For $Q$-partitions $A,B$ of $X$, let $A\leq_WB$ mean that there is a  continuous function $f$ on $X$ such that $A(x)\leq_QB(f(x))$ for each $x\in X$. The case of sets corresponds to the case of 2-partitions. Let  ${\mathbf B}(Q^X)$ be the set of Borel $Q$-partitions $A$ (for which  $A^{-1}(q)\in{\mathbf B}(X)$ for all $q\in Q$). A celebrated theorem of van Engelen, Miller and Steel (see Theorem 3.2 in \cite{ems87}) shows that if $Q$ is a countable better quasiorder (bqo) then  $\mathcal{W}_Q=({\mathbf B}(Q^\mathcal{N});\leq_W)$  is a bqo. Although this theorem gives an important information about the quotient-poset of $\mathcal{W}_Q$, it is  far from a  characterisation.

Many efforts (see e.g. \cite{he93,s07a,s16,s17} and references therein) to characterise the quotient-poset of $\mathcal{W}_Q$ were devoted to {\em $k$-partitions of $\mathcal{N}$}.    Our approach in \cite{s07a,s16,s17} to this problem was to characterise the initial segments $({\mathbf \Delta}^0_\alpha(k^\mathcal{N});\leq_W)$ for bigger and bigger ordinals $2\leq\alpha<\omega_1$. To achieve this, we defined  structures of iterated labeled trees and forests with the so called homomorphism quasiorder and discovered  useful properties of some natural operations on the iterated labeled  forests and on $Q$-partitions.

An important progress was recently achieved in \cite{km17} where a full characterisation of the quotient-poset of $\mathcal{W}_Q$ for arbitrary countable bqo $Q$ is obtained, using an extended set of iterated labeled trees $(\mathcal{T}_{\omega_1}(Q);\leq_h)$ with the homomorphism quasiorder $\leq_h$. Namely, $(\mathcal{T}_{\omega_1}(Q);\leq_h)$ is equivalent to the substructure of $\mathcal{W}_Q$ formed by the $\sigma$-join-irreducible elements   (the equivalence means isomorphism of the corresponding quotient-posets) via an embedding $\mu:\mathcal{T}_{\omega_1}(Q)\to\mathcal{W}_Q$. 
The Wadge hierarchy of $Q$-partitions of $\calN$ may be thus written as the family $\{\mathcal{W}_Q(T)\}_{T\in\mathcal{T}_{\omega_1}(Q)}$, where $\mathcal{W}_Q(T)=\{A\in Q^\mathcal{N}\mid A\leq_W\mu(T)\}$, and it exhausts all principal ideals of $\mathcal{W}_Q$ formed by $\sigma$-join-irreducible $Q$-Wadge degrees. For the case of 2-partitions this yields a new characterization of the Wadge hierarchy of sets.

Our definition of the $Q$-IFH may be now  sketched as follows. In arbitrary space $X$ (and even in a more general situation) we define the family $\{\mathcal{L}(X,T)\}_{T\in\mathcal{T}_{\omega_1}(Q)}$ of classes of $Q$-partitions of $X$ which we call the $Q$-IFH  in $X$. We then show that if $X$ is quasi-Polish then $\mathcal{L}(X,T)=\{A:X\to Q\mid A\circ\xi\in\mathcal{W}_Q(T)\}$ for all $T\in\mathcal{T}_{\omega_1}(Q)$ (at least for $Q=\bar{k}$). For the case of $2$-partitions we obtain a set-theoretic characterisation of the Wadge hierarchy of sets defined above within the second approach. This characterisation looks rather different from a set-theoretic description of the Wadge hierarchy in \cite{wad84} (see also \cite{lo12}). Note that the characterisations in \cite{wad84,lo12} cannot be straightforwardly extended to arbitrary spaces since they use specific features of the Baire space. The properties of $Q$-IFH in $X$ strongly depend on $Q$ (we distinguish the cases when $Q$ is an arbitrary quasiorder, a bqo, an antichain, $Q=\bar{k}$, $Q=\bar{2}$) and on $X$ (we distinguish the cases when $X$ is a set, an arbitrary space, a second countable space, a quasi-Polish space, the Baire space), which is reflected in many formulations below. 

 Having  papers \cite{s12,s16,s19} at hand would probably simplify reading of the present paper because they contain simpler versions of some notions and results  based on similar ideas. The main technical notions for the infinitary case are a bit more complicated than for the finitary case (considered e.g. in \cite{s12,s19}) but the ideas are the same. 

After recalling necessary preliminaries in the next section, we define  in Section \ref{qfine} the $Q$-IFH  and establish its general properties. In Section \ref{char} we prove additional properties of the $Q$-IFH in second countable spaces and in quasi-Polish spaces. In particular, we prove the above-mentioned preservation property and Hausdorff-Kuratowski-type theorems and show that in the Baire space the $Q$-IFH coincides with the $Q$-Wadge hierarchy from \cite{km17}. We also examine when levels of this hierarchy have natural representations, are downward closed under Wadge reducibility and have Wadge complete $Q$-partition. 

%In particular, the proof of a broad extension of the Steel separation theorem becomes a straightforward set-theoretic exercise instead of a sophisticated use of  specially designed Gale-Stewart games. 

In Section \ref{effect} we also briefly discuss the effective finitary version of Wadge hierarchy developed in \cite{s19} and its relation to the non-effective version developed here. We conclude in Section \ref{con} with some of the remaining open questions.

\section{Preliminaries}\label{prel}

In this section we briefly recall some notation, notions and facts used throughout the paper. Some more special information is recalled in the corresponding sections below.

\subsection{Well and better quasiorders}\label{bqo}

We use standard set-theoretical notation. In particular, $Y^X$ is the set of functions from $X$ to $Y$,
$P(X)$ is the class of subsets of a set $X$, $\check{\mathcal{C}}$ is the class of complements $X\setminus C$ of sets $C$ in $\mathcal{C}\subseteq P(X)$.
We assume the reader to be acquainted with the notion of ordinal (see e.g.  \cite{km67}). Ordinals are denoted by $\alpha, \beta, \gamma,\ldots$.  Every ordinal $\alpha$ is the set
of  smaller ordinals, in particular $ k=\{0,1,\ldots,k-1\}$ for each $k<\omega$, and $ \omega=\{0,1,2,\ldots\}.$ 
We use some notions and facts of ordinal arithmetic. In particular,
$\alpha+\beta$, $\alpha\cdot\beta$ and $\alpha^{\beta}$ denote the
ordinal addition, multiplication and exponentiation of $\alpha$ and
$\beta$, respectively. Every positive ordinal $\alpha$ is uniquely representable in the form $\alpha=\omega^{\alpha_0}+\cdots+\omega^{\alpha_n}$ where $n<\omega$ and $\alpha\geq\alpha_0\geq\cdots\geq\alpha_n$; we denote $\alpha*=\omega^{\alpha_0}$.
The first
non-countable ordinal is denoted by $\omega_1$.  

We use  standard notation and terminology on partially ordered sets (posets).  Recall that a
{\em quasiorder} (qo) is a structure $(P;\leq)$ satisfying the axioms of
reflexivity $\forall x(x\leq x)$ and transitivity $\forall x\forall
y\forall z(x\leq y\wedge y\leq z\to x\leq z)$. Any qo $\leq$ on $P$ induces the equivalence relation defined by $a\equiv b\leftrightarrow a\leq
b\wedge b\leq a$. The corresponding quotient structure of $(P;\leq)$ is called  {\em the
quotient-poset} of $P$. To avoid complex notation, we sometimes abuse terminology about posets by applying it also to qo's; in
such cases we just mean the corresponding quotient-poset.

A qo $(P;\leq)$ is  {\em well-founded} if it has no
infinite descending chains $a_0>a_1>\cdots$. In this case there are a unique ordinal
$rk(P)$ and a unique rank function $rk_P$ from $P$ onto $rk(P)$
satisfying $a<b\to rk(a)<rk(b)$. It is defined by induction
$rk_P(x)=sup\{rk_P(y)+1\mid y<x\}. $ The ordinal $rk(P)$ is called
the {\em rank} (or {\em height}) of $P$, and the ordinal $rk_P(x)$ is called the {\em rank of  $x\in P$ in $P$}.   

A {\em well quasiorder} (wqo) is a qo $Q=(Q;\leq_Q)$ that has neither infinite
descending chains nor infinite antichains. Although  wqo's are closed under many natural finitary constructions like forming finite labeled words or trees, they are not always closed under important infinitary constructions.
Nevertheless, it turns out possible to find a natural subclass of wqo's, called better quasiorders (bqo's) which contains most of the ``natural'' wqo's (in particular, all finite qo's) and has strong closure properties also for many infinitary constructions.  The notion of bqo is due to C. Nash-Williams. We omit a bit technical
notion of bqo which is used only in formulations. For more details on bqo's, we refer the reader to \cite{si85}. 

Recall that
{\em semilattice} is a structure
$(S;\sqcup)$ with binary operation $\sqcup$ such that $(x\sqcup
y)\sqcup z=x\sqcup (y\sqcup z)$, $x\sqcup y= y\sqcup x$  and $x\sqcup x=x$, for
all $x,y,z\in S$. By $\leq$ we denote the induced partial order on $S$:
$x\leq y$ iff $x\sqcup y=y$. The operation $\sqcup$ can be recovered from $\leq$ since $x\sqcup y$ is the supremum of $x,y$ w.r.t. $\leq$.
By {\em  $\sigma$-semilattice} we mean a semilattice where also  supremums $\bigsqcup y_j=y_0\sqcup y_1\sqcup\cdots$ of countable sequences of elements $y_0,y_1,\ldots$ exist. Element $x$ of a $\sigma$-semilattice $S$ is
{\em $\sigma$-join-irreducible} if it cannot be represented as the countable supremum of elements strictly below $x$. As first stressed in \cite{s07}, the $\sigma$-join-irreducible elements play a central role in the study of Wadge degrees of $k$-partitions. The same applies to several variations of Wadge degrees, including the Wadge degrees of $Q$-partitions for a countable bqo $Q$.

\subsection{Classical hierarchies in topological spaces}\label{chier}

We assume the reader to be familiar with  basic notions of
topology \cite{en89}. The underlying set  of a topological space $X$  will  be usually also denoted by $X$, in abuse
of notation. We usually abbreviate ``topological space'' to
``space''. A space is \emph{zero-dimensional} if it has a basis of
clopen sets. Recall that a \emph{basis} for the topology on $X$ is a
set $\cal B$ of open subsets of $X$ such that for every $x\in X$ and
open $U$ containing $x$ there is $B\in \cal B$ satisfying $x\in
B\subseteq U$. We sometimes shorten ``countably based $T_0$-space'' to ``cb$_0$-space''.

Let $\omega$ be the space of non-negative integers with the
discrete topology. 
Let $\calN=\omega^\omega$ be the set of all infinite
sequences of natural numbers (i.e., of all functions $x \colon
\omega \to \omega$). Let $\omega^*$ be the set of finite sequences
of elements of $\omega$, including the empty sequence $\varepsilon$. For
$\sigma\in\omega^*$ and $x\in\calN$, we write
$\sigma\sqsubseteq x$ to denote that $\sigma$ is an initial
segment of the sequence $x$. By $\sigma x=\sigma\cdot x$ we
denote the concatenation of $\sigma$ and $x$, and by
$\sigma\cdot\calN$ the set of all extensions of $\sigma$ in
$\calN$. For $x\in\calN$, we can write
$x=x(0)x(1)\dotsc$ where $x(i)\in\omega$ for each $i<\omega$. For
$x\in\calN$ and $n<\omega$, let $x\upharpoonright n=x(0)\dotsc x(n-1)$
denote the initial segment of $x$ of length $n$. By endowing $\calN$ with the product of the discrete
topologies on $\omega$, we obtain the so-called \emph{Baire space}.
The product topology coincides with the topology
generated by the collection of sets of the form
$\sigma\cdot\calN$ for $\sigma\in\omega^*$. It is well known that $\calN\times\calN$ and $\calN^\omega$ are homeomorphic to $\calN$.

A {\em tree} is a non-empty set $T\subseteq\omega^*$ which is closed downwards under the prefix relation $\sqsubseteq$. The empty string $\varepsilon$ is the {\em root} of any tree. A {\em leaf} of $T$ is a maximal element of $(T;\sqsubseteq)$.  A tree is {\em pruned} if it has no leafs. A {\em path through} a tree $T$ is an element $x\in\calN$ such that $x\upharpoonright n\in T$ for each $n\in\omega$. For any tree and any $\tau\in T$, let $[T]$ be the set of paths through $T$ and $T(\tau)=\{\sigma\mid\tau\sigma\in T\}$. The non-empty closed subsets of $\mathcal{N}$ coincide with the sets $[T]$ where $T$ is pruned; every nonempty closed set is a retract of $\mathcal{N}$. 

We call a tree $T$ {\em normal} if $\tau(i+1)\in T$ imply $\tau i\in T$. A tree  is {\em infinite-branching} if with every non-leaf node $\tau$ it contains all its successors $\tau i$; every infinite branching tree is normal.
A tree is {\em well founded} if there is no path through it (i.e., $(T;\sqsupseteq)$ is well founded). The rank of the latter poset is called the rank of $T$; the ranks of well founded trees are precisely the countable ordinals. By a {\em forest} we mean a set of strings $T\setminus\{\varepsilon\}$, for some tree $T$; usually we assume forests to be non-empty. Sometimes we use  other obvious notation on trees. E.g. with any sequence of trees $\{T_0,T_1,\ldots\}$ we associate the tree $T=\{\varepsilon\}\cup0\cdot T_0\cup1\cdot T_1\cup\cdots$ such that $T(i)=T_i$ for each $i<\omega$.

A \emph{pointclass} in a space $ X $ is  a class $\bfGamma(X)\subseteq P(X) $ of subsets of $ X $; let $\check{\bfGamma}(X)=\{A\subseteq X\mid X\setminus A\in\bfGamma(X)\}$. 
A \emph{family of pointclasses} \cite{s13} is a family $ \bfGamma=\{\bfGamma(X)\}_X $ 
indexed by arbitrary topological spaces $X$ (or by spaces in a reasonable class) such that each $ \bfGamma(X) $ is
a pointclass in $ X $ and $ \bfGamma $ is closed under
continuous preimages, i.e. $ f^{-1}(A)\in\bfGamma(X) $
for every $ A\in\bfGamma(Y) $ and every continuous function $ f \colon X\to Y $. 
A basic example of a family of pointclasses is given by the family
$\mathcal{O}=\{\tau_X\}_X$ of  topologies in arbitrary spaces $X$.

We will use the following operations on families of pointclasses: 
the operation $\bfGamma\mapsto\bfGamma_\sigma$,
where $\bfGamma(X)_\sigma$ is the set of all countable unions of sets in $\bfGamma(X)$, 
the operation $\bfGamma\mapsto\bfGamma_\delta$,
where $\bfGamma(X)_\delta$ is the set of all countable intersections of sets in $\bfGamma(X)$, 
the operation $\bfGamma\mapsto\bfGamma_c$,
where $\bfGamma(X)_c=\check{\bfGamma}(X)$, 
the operation $\bfGamma\mapsto\bfGamma_d$,
where $\bfGamma(X)_d$ is the set of all differences of sets in $\bfGamma(X)$, 
the operation $\bfGamma\mapsto\bfGamma_\exists$ defined by 
$\bfGamma_\exists(X):=\{\exists^\mathcal{N}(A)\mid A\in\bfGamma(\calN\times X)\}$,
where $\exists^\mathcal{N}(A):=\{x\in X\mid \exists p\in\calN.(p,x)\in A\}$ is the projection of $A\subseteq\calN\times X$ 
along the axis $\calN$,
and, finally, the operation $\bfGamma\mapsto\bfGamma_\forall$ defined by 
$\bfGamma_\forall(X):=\{ \forall^\mathcal{N}(A) \mid A\in\bfGamma(\calN\times X)\}$,
where $\forall^\mathcal{N}(A):=\{x\in X\mid \forall p\in\calN.(p,x)\in A\}$.

The  operations on families of pointclasses enable to provide short uniform descriptions 
of the classical hierarchies in arbitrary spaces. 
E.g., the {\em Borel  hierarchy} is the sequence of families of pointclasses
$\{\bfSig^0_\alpha\}_{ \alpha<\omega_1}$ defined by induction on $\alpha$ as follows \cite{s06,br}:
$\bfSig^0_0(X):=\{\emptyset\}$,
$\bfSig^0_1 := \mathcal{O}$ (the family of open sets), $\bfSig^0_2 := (\bfSig^0_1)_{d\sigma}$,
and  $\bfSig^0_\alpha(X)
=(\bigcup_{\beta<\alpha}\bfSig^0_\beta(X))_{c\sigma}$ for
$\alpha>2$.
The sequence $\{\bfSig^0_\alpha(X)\}_{ \alpha<\omega_1}$ is called \emph{the Borel hierarchy} in $X$.
We also set $\bfPi^0_\beta(X) = (\bfSig^0_\beta(X))_c $
and $\bfDelta^0_\alpha(X) = \bfSig^0_\alpha(X) \cap \bfPi^0_\alpha (X)$.
The classes $\bfSig^0_\alpha(X),\bfPi^0_\alpha(X),{\bfDelta}^0_\alpha(X)$
are called  \emph{levels} of the Borel hierarchy in $X$. 
 The class $\mathbf{B}(X)$ of \emph{Borel sets} in $X$ is defined as the union of all levels 
 of the Borel hierarchy in $X$; it coincides with the smallest $\sigma$-algebra of subsets of $X$ containing the open sets.
We have $\bfSig^0_\alpha(X)\cup\bfPi^0_\alpha(X)\subseteq\bfDelta^0_\beta(X)$ for all $\alpha<\beta<\omega_1$. 

 The {\em hyperprojective hierarchy} is the sequence of families of pointclasses
 $\{\bfSig^1_\alpha\}_{ \alpha<\omega_1}$ defined by induction on $\alpha$ as follows:  
 $\bfSig^1_0=\bfSig^0_2$, 
 $\bfSig^1_{\alpha+1}=(\bfSig^1_\alpha)_{c\exists}$, 
 $\bfSig^1_\lambda=(\bfSig^1_{<\lambda})_{\delta \exists}$, 
 where $\alpha,\lambda<\omega_1$, $\lambda$ is a limit ordinal, 
 and $\bfSig^1_{<\lambda}(X)=\bigcup_{\alpha<\lambda}\bfSig^1_\alpha(X)$.
 In this way, we obtain for any cb$_0$-space $X$ the sequence 
 $\{\bfSig^1_\alpha(X)\}_{\alpha<\omega_1}$, which we call here \emph{the hyperprojective hierarchy in $X$}. 
 The pointclasses $\bfSig^1_\alpha(X)$, $\bfPi^1_\alpha(X)=(\bfSig^1_\alpha(X))_c$ and 
 $\bfDelta^1_\alpha(X)=\bfSig^1_\alpha(X)\cap\bfPi^1_\alpha(X)$ are called 
 \emph{levels of the hyperprojective hierarchy in $X$}. 
 The finite non-zero levels of the hyperprojective hierarchy coincide with the corresponding levels 
 of the Luzin projective hierarchy. 
 
 We do not recall the well known definition of the Hausdorff difference hierarchy over $\bfSig^0_{\alpha}(X)$, $\alpha\geq1$, which is denoted by $\{D_\beta(\bfSig^0_{\alpha}(X))\}_{\beta<\omega_1}$ or by $\{\bfSig^{-1,\alpha}_\beta(X)\}_{\beta<\omega_1}$. The definitions may be found e.g. in \cite{ke95,s16}.

We recall some structural properties of pointclasses  (see e.g. \cite{ke95}). 

\begin{Definition}\label{d-in-po3}
\begin{enumerate}\itemsep-1mm

\item  A pointclass $\bfGamma(X)$ in $X$ has the {\em separation property}  if for
every two disjoint sets $A,B\in\bfGamma(X)$ there is a set
$C\in\bfGamma(X)\cap \check{\bfGamma}(X)$ with $A\subseteq
C\subseteq X\setminus B$. 

\item  A pointclass $\bfGamma(X)$ has the {\em reduction property} i.e. for
all $C_0, C_1\in\bfGamma(X)$ there are disjoint
$C^\prime_0,C^\prime_1\in\bfGamma(X)$ such that
$C^\prime_i\subseteq C_i$ for  $i<2$ and $C_0\cup
C_1=C^\prime_0\cup C^\prime_1$. The pair $(C^\prime_0,C^\prime_1)$
is called a reduct for the pair $(C_0,C_1)$.

\item  A pointclass $\bfGamma(X)$ in $X$ has the {\em $\sigma$-reduction property}
if for each countable sequence $C_0, C_1,\ldots$ in ${\mathcal A}$
there is a countable sequence $C^\prime_0, C^\prime_1,\ldots$ in
$\bfGamma(X)$ (called a reduct of $C_0, C_1,\ldots$) such that
$C^\prime_i\cap C^\prime_j=\emptyset$ for all $i\not=j$ and
$\bigcup_{i<\omega}C^\prime_i=\bigcup_{i<\omega}C_i$.

\end{enumerate}

\end{Definition}

It is well-known that if $\bfGamma(X)$ has the
reduction property then the dual class $\check{\bfGamma}(X)$ has the
separation property, but not vice versa, and that
if $\bfGamma(X)$ has the $\sigma$-reduction property then
$\bfGamma(X)$ has the reduction property but not vice versa.
Let $X$ be a cb$_0$-space. It is known (see e.g. \cite{ke95,s13}) that any level $\bfSig^0_{2+\alpha}(X)$, $\alpha<\omega_1$, has the $\sigma$-reduction property, and if $X$ is zero-dimensional then also $\bfSig^0_{1}(X)$ has the $\sigma$-reduction property.

\subsection{Quasi-Polish spaces and admissible representations}\label{qpolish}

A space $X$ is \emph{Polish} if it is countably based
and metrizable with a metric $d$ such that $(X,d)$ is a complete
metric space. Examples of Polish spaces are $\omega$,
$\calN$, the Cantor space $\mathcal{C}$, the space of reals $\IR$ and its Cartesian
powers $\IR^n $ ($ n < \omega $), the closed unit interval $ [0,1]
$, the Hilbert cube $ [0,1]^\omega $ and the space
$\IR^\omega $. 

Quasi-Polish spaces were identified and thoroughly studied by M. de Brecht \cite{br} (see also \cite{chen} for additional information). Informally, this is a natural class of spaces which contains all Polish spaces, many important non-Hausdorff spaces (like  $\omega$-continuous domains) and has essentially the same DST as Polish spaces.
Let $\Pomega$ be the space of subsets of $\omega$ equipped
with the Scott topology, a countable basis of which is
formed by the sets $\{A\subseteq\omega\mid F\subseteq A\}$, where
$F$ ranges over the finite subsets of $\omega$. By a {\em quasi-Polish space} we mean a space homeomorphic to a $\bfPi^0_2$-subspace of $\Pomega$.  There are several interesting characterizations of quasi-Polish spaces. For this paper the following characterization in terms of representations is relevant.

A \emph{representation} of a set $X$ is a surjection from a subspace of  $\calN$ onto $X$. Such a representation is {\em total} if its domain is $\calN$. Representation $\mu$ is {\em (continuously) reducible} to a representation $\nu$ ($\mu\leq_c\nu$) if $\mu=\nu\circ f$ for some continuous partial function $f$ on $\calN$. Representations $\mu,\nu$ are {\em (continuously) equivalent} ($\mu\equiv_c\nu$) if $\mu\leq_c\nu$ and $\nu\leq_c\mu$.  
A basic notion of Computable
Analysis \cite{wei00} is the notion of admissible representation. A
representation $\mu$ of a space $X$ is \emph{admissible}, if it is
continuous and any continuous function $\nu:Z \to X$ from a subset
$Z\subseteq\calN$ to $X$ is continuously reducible to $\mu$. Clearly, any two admissible representations of a space are continuously equivalent. As shown in \cite{bh02}, any continuous open surjection from a subspace of $\calN$ onto $X$ is an admissible representation of $X$.
In \cite{br} the following  characterization of
quasi-Polish spaces  was obtained:

\begin{Proposition}\cite{br}\label{adm}
A cb$_0$-space $X$ is quasi-Polish iff it has a total admissible representation iff there is a continuous open surjection from $\calN$ onto $X$.
\end{Proposition}

The classical Borel, Luzin and Hausdorff hirarchies in quasi-Polish spaces have properties very similar to their properties in Polish spaces \cite{br}. In particular,
for any uncountable quasi-Polish space $X$ and any $\alpha<\omega_1$, $\bfSig^0_\alpha(X)\not\subseteq\bfPi^0_\alpha(X)$ and $\bfSig^1_\alpha(X)\not\subseteq\bfPi^1_\alpha(X)$. For any quasi-Polish space $X$, the Suslin theorem $\bigcup_{\alpha<\omega_1}\bfSig^0_{1+\alpha}(X)=\mathbf{B}(X)=\bfDelta^1_1(X)$ and the Hausdorff-Kuratowski theorem \cite{ke95,br} (saying that $\bigcup_{\beta<\omega_1}\bfSig^{-1,\alpha}_\beta(X)=\bfDelta^0_{\alpha+1}(X)$ for all $\alpha\geq1$) hold.

Quasi-Polish spaces also share  properties of Polish spaces related to Baire category \cite{ke95,chen}. According to Corollary 52 in \cite{br} (see also \cite{chen}), every quasi-Polish space $X$ is completely Baire, in particular every nonempty closed set $F\subseteq X$ is non-meager in $F$. Using the technique of category quantifiers \cite{ke95}, one can show the following preservation property \cite{sr07,br} of levels of the Borel hierarchy. 

\begin{Proposition}\cite{sr07,br}\label{preserve}
Let $f:X\to Y$ be a continuous open surjection between cb$_0$-spaces, $\alpha<\omega_1$, and $A\subseteq Y$. Then $A\in\bfSig^0_{1+\alpha}(Y)$ iff $f^{-1}(A)\in\bfSig^0_{1+\alpha}(X)$. Also, every fiber $f^{-1}(y)$ is quasi-Polish, hence non-meager in $f^{-1}(y)$.
\end{Proposition}

%A similar preservation property holds for the levels of Hausdorff-Kuratowski difference hierarchies \cite{br}. 

\subsection{Wadge hierarchy in $\mathcal{N}$}\label{wadgedst}

Here we give some additional information on the Wadge hierarchy in the Baire space mentioned in the Introduction.  
In \cite{wad84} W. Wadge  (using the Borel
determinacy) proved the following result:
The structure $({\mathbf B}(\calN);\leq_W)$ of  Borel sets in the Baire
space is semi-well-ordered (i.e., it is well-founded and for all
$A,B\in{\mathbf B}(\mathcal{N})$ we have $A\leq_WB$ or $\overline{B}\leq_WA$). In particular, there is no antichain of size 3 in $({\mathbf B}(\calN);\leq_W)$. He has also computed the rank $\upsilon$ of $({\mathbf B}(\calN);\leq_W)$ which we call the Wadge ordinal. Recall that a set $A$ is {\em self-dual} if $A\leq_W\overline{A}$. W. Wadge has shown that
if a Borel set  is self-dual (resp. non-self-dual) then any Borel set of the next Wadge rank is non-self-dual (resp. self-dual), a Borel set of Wadge rank of countable cofinality is self-dual, and a Borel set of Wadge rank of uncountable cofinality is non-self-dual. This characterizes the  structure of Wadge degrees of Borel sets up to isomorphism.

In \cite{ste80} the following separation theorem for the Wadge hierarchy was established:
For any non-self-dual Borel set $A$   exactly one of the principal ideals
$\{X\mid X\leq_WA\}$, $\{X\mid X\leq_W\overline{A}\}$ has the
separation property. 

The mentioned results give rise to the {\em Wadge hierarchy} which
is, by definition, the sequence
$\{\bfSig_\alpha(\calN)\}_{\alpha<\upsilon}$ (where $\upsilon$ is the Wadge ordinal) of all non-self-dual principal
ideals of $({\bf B}(\calN);\leq_W)$ that do not have the separation
property  and satisfy for all $\alpha<\beta<\upsilon$ the strict
inclusion $\bfSig_\alpha(\calN)\subset\bfDelta_\beta(\calN)$ where, as usual, $\bfDelta_\alpha(\calN)=\bfSig_\alpha(\calN)
\cap\bfPi_\alpha(\calN)$. 

The Wadge hierarchy subsumes the classical hierarchies in the Baire space, in particular
${\bfSig}_\alpha(\calN)={\bfSig}^{-1}_\alpha(\calN)$ for each
$\alpha<\omega_1$, 
${\bfSig}_1(\calN)={\bfSig}^0_1(\calN)$,
${\bfSig}_{\omega_1}(\calN)={\bfSig}^0_2(\calN)$,
${\bfSig}_{\omega_1^{\omega_1}}(\calN)={\bfSig}^0_3(\calN)$ and so on.
Thus, the sets  of finite Borel rank coincide with the sets of Wadge
rank less than
$\lambda=sup\{\omega_1,\omega_1^{\omega_1},\omega_1^{(\omega_1^{\omega_1})},\ldots\}$.
Note that $\lambda$ is the smallest solution of the ordinal equation
$\omega_1^\varkappa=\varkappa$. Hence, the reader should carefully distinguish $\bfSig_\alpha(\calN)$ and $\bfSig^0_\alpha(\calN)$. To give the reader an impression about the Wadge ordinal we note that
the rank of the qo $(\bfDelta^0_\omega(\calN);\leq_W)$ is the
$\omega_1$-st  solution of the ordinal equation
$\omega_1^\varkappa=\varkappa$ \cite{wad84}.

We summarise some properties of the Wadge hierarchy of sets in the Baire space which will be tested for survival under generalisations to cb$_0$-spaces (or to quasi-Polish spaces) and to $Q$-partitions below:

\begin{enumerate}
 \item The levels of the Wadge hierarchy are semi-well-ordered by inclusion.
 \item The Wadge hierarchy does not collapse, i.e. $\bfSig_\alpha\not\subseteq\bfPi_\alpha$ for all $\alpha<\upsilon$. 
 \item The Wadge degrees of Borel sets coincide with the sets $\bfSig_\alpha\setminus\bfPi_\alpha$, $\bfPi_\alpha\setminus\bfSig_\alpha$, $\bfDelta_{\alpha+1}\setminus(\bfSig_\alpha\cup\bfPi_\alpha)$ (where $\alpha<\upsilon$), and $\bfDelta_{\lambda}\setminus (\bigcup_{\alpha<\lambda}\bfSig_\alpha)$ (where $\lambda<\upsilon$ is a limit ordinal of countable cofinality).
 \item If $\lambda<\upsilon$ is a limit ordinal of uncountable cofinality then $\bfDelta_{\lambda}= \bigcup_{\alpha<\lambda}\bfSig_\alpha$.
 \item All levels are downward closed under Wadge reducibility.
 \item The levels in item (3) are precisely those having  Wadge complete sets.
% \item The levels $\bfPi_\alpha$ (but not the  levels $\bfSig_\alpha$) have the separation property.
\end{enumerate}

\section{Infinitary fine hierarchies in a set}\label{qfine}

In this section we define the infinitary fine hierarchy  and prove some of its basic properties. The $Q$-partition version of this hierarchy will be called the $Q$-IFH, for abbreviation. Definitions and results in this section extend (and in fact simplify) the corresponding material from Section 5 in \cite{s16}. Three following subsections describe some related technical notions.

\subsection{Iterated $Q$-trees}\label{trees}

Here we describe a notation system for  levels of the $Q$-IFH.

For any qo $Q$, a {\em $Q$-tree} is a pair $(T,t)$ consisting of an infinite-branching well founded tree $T\subseteq\omega^*$ and a labeling $t:T\to Q$. Let $\mathcal{T}(Q)$ be the set of $Q$-trees quasi-ordered by the relation: $(T,t)\leq_h(V,v)$ iff there is a monotone function $\varphi:T\to V$ with $\forall v\in T(t(x)\leq_Qv(\varphi(x)))$; such a function $\varphi$ is called a {\em morphism from $(T,t)$ to $(V,v)$}. As follows from a Laver's theorem on bqo's, if $Q$ is bqo then so is also $(\mathcal{T}(Q);\leq_h)$ which is usually shortened to $\mathcal{T}(Q)$. Thus, $\mathcal{T}$ is an operator on the class BQO of all bqo's. 
The operator $\mathcal{T}$  and its iterates like $\mathcal{T}\circ\mathcal{T}\circ\mathcal{T}$ were introduced in \cite{s07,s16} and turned out crucial for characterising some initial segments of $\mathcal{W}_{\bar{k}}$  \cite{s16,s17}. 

In \cite{km17} a more powerful iteration procedure  was invented which yields the set $\mathcal{T}_{\omega_1}(Q)$ of labeled trees sufficient for characterising $\mathcal{W}_Q$, as discussed in the Introduction. We give a slightly different (but equivalent) definition of $\mathcal{T}_{\omega_1}(Q)$ which is more convenient for our purposes here. The differences are caused by our desire to first work only with trees (introducing forests at the last stage, see below), and to relate the qo $\unlhd$ (defined below) to the qo $\leq_h$.

Let $\sigma=\sigma(Q,\omega_1)=\{q,s_\alpha,F_q,F_\alpha\mid q\in Q,\alpha<\omega_1\}$ be the signature consisting of constant symbols $q$, unary function symbols $s_\alpha$, and $\omega$-ary function symbols $F_q,F_\alpha$ (of course we assume all signature symbols to be distinct, in particular $Q\cap\omega_1=\emptyset$). Let $\mathbb{T}_\sigma$ be the set of $\sigma$-terms without variables obtained by the standard inductive definition: Any constant symbol $q$ is a term; if $u$ is a term then so is also $s_\alpha(u)$; if $u_0,u_1,\ldots$ are terms then so are also $F_q(u_0,\ldots), F_{\alpha}(u_0,\ldots)$. Informally, $F_q(u_0,\ldots)$ and $F_{\alpha}(u_0,\ldots)$ are interpreted as $q\to(u_0\sqcup\cdots)$ and $s_\alpha(u_0)\to(u_1\sqcup\cdots)$ respectively (cf. \cite{km17} where e.g. the first expression denotes the tree $\varepsilon\cup0\cdot u_0\cup\cdots$ with the root $\varepsilon$  labeled by $q$), hence our modification simply avoids forests from the inductive definition. 

The $\sigma$-terms are represented by (or even identified with) the normal well founded trees with constants on the leafs and other signature symbols  on the non-leaf nodes such that the nodes labeled with $s_\alpha$ have  the unique successor while the nodes labeled by $F_q$ of $F_\alpha$ have all successors. Such {\em syntactic trees} enable definitions and proofs by induction on terms (i.e., on the  ranks of syntactic trees) because the subterms of (the syntactic tree of) $u$ are precisely the trees $u(\tau)$, see Section \ref{chier}. Obviously, the set $\mathbb{T}_\sigma$ is partitioned into three parts: {\em constant terms} (i.e., the terms $q$ for some $q\in Q$), {\em $s$-terms} (i.e., the terms $s_\alpha(u)$ for unique $\alpha<\omega_1$ and $u\in \mathbb{T}_\sigma$) and {\em $F$-terms} (i.e., the terms $F_q(u_0,\ldots)$ or $F_\alpha(u_0,\ldots)$ for unique $q\in Q$, $\alpha<\omega_1$, and $u_0,u_1,\ldots\in\mathbb{T}_\sigma$). 
We define by induction on terms the binary relation $\unlhd$ on $\mathbb{T}_\sigma$ as follows (cf. Definition 3.1 and its extensions in \cite{km17}). The relation $\unlhd$ on $\mathbb{T}_\sigma$ is in fact equivalent to the relation $\unlhd$ in \cite{km17} restricted to the tree-terms. 

\begin{Definition}\label{treeterms}
\begin{enumerate}
 \item $q\trianglelefteq r$ iff $q\leq_Qr$;
 \item $q\trianglelefteq s_\alpha(u)$ iff $q\unlhd u$;
 \item $q\trianglelefteq F_r(u_0,\ldots)$ iff $q\trianglelefteq r$ or $q\trianglelefteq u_i$ for some $i\geq0$;
 \item $q\trianglelefteq F_\alpha(u_0,\ldots)$ iff  $q\trianglelefteq u_i$ for some $i\geq0$;
 \item $s_\alpha(u)\trianglelefteq r$ iff $u\trianglelefteq r$;
 \item $s_\alpha(u)\trianglelefteq s_\beta(v)$ iff ($\alpha<\beta$ and $u\trianglelefteq s_\beta(v)$) or ($\alpha=\beta$ and $u\trianglelefteq v$) or ($\alpha>\beta$ and $s_\alpha(u)\trianglelefteq v$);
 \item $s_\alpha(u)\trianglelefteq F_r(v_0,\ldots)$ iff $s_\alpha(u)\trianglelefteq r$ or $s_\alpha(u)\trianglelefteq v_i$ for some $i\geq0$;
 \item $s_\alpha(u)\trianglelefteq F_\beta(v_0,\ldots)$ iff $s_\alpha(u)\trianglelefteq s_\beta(v_0)$ or $s_\alpha(u)\trianglelefteq v_i$ for some $i\geq1$;
 \item $F_q(u_0,\ldots)\trianglelefteq r$ iff $q\trianglelefteq r$ and $u_i\trianglelefteq r$ for all $i\geq0$;
 \item $F_q(u_0,\ldots)\trianglelefteq s_\alpha(v)$ iff $q\trianglelefteq s_\alpha(v)$ and $u_i\trianglelefteq s_\alpha(v)$ for all $i\geq0$;
 \item $F_q(u_0,\ldots)\trianglelefteq F_r(v_0,\ldots)$ iff ($q\trianglelefteq r$ and $u_i\trianglelefteq F_r(v_0,\ldots)$ for all $i\geq0$) or $F_q(u_0,\ldots)\trianglelefteq v_i$ for some $i\geq0$; 
 \item $F_q(u_0,\ldots)\trianglelefteq F_\beta(v_0,\ldots)$ iff ($q\trianglelefteq s_\beta(v_0)$ and $u_i\trianglelefteq F_\beta(v_0,\ldots)$ for all $i\geq0$) or $F_p(u_0,\ldots)\trianglelefteq v_i$ for some $i\geq1$;
  \item $F_\alpha(u_0,\ldots)\trianglelefteq r$ iff $s_\alpha(u_0)\trianglelefteq r$ and $u_i\trianglelefteq q$ for all $i\geq1$;
 \item $F_\alpha(u_0,\ldots)\trianglelefteq s_\beta(v)$ iff $s_\alpha(u_0)\trianglelefteq s_\beta(v)$ and $u_i\trianglelefteq s_\beta(v)$ for all $i\geq1$;
 \item $F_\alpha(u_0,\ldots)\trianglelefteq F_r(v_0,\ldots)$ iff ($s_\alpha(u_0)\trianglelefteq r$ and $u_i\trianglelefteq F_q(v_0,\ldots)$ for all $i\geq1$) or $F_\alpha(u_0,\ldots)\trianglelefteq v_i$ for some $i\geq0$; 
 \item $F_\alpha(u_0,\ldots)\trianglelefteq F_\beta(v_0,\ldots)$ iff ($s_\alpha(u_0)\trianglelefteq s_\beta(v_0)$ and $u_i\trianglelefteq F_\beta(v_0,\ldots)$ for all $i\geq1$) or $F_\alpha(u_0,\ldots)\trianglelefteq v_i$ for some $i\geq1$.
\end{enumerate}
\end{Definition}

The remarks above and the arguments in \cite{km17} show that $(\mathbb{T}_\sigma;\trianglelefteq)$ is a bqo. Let $\mathbb{T}_{q,s}$ be the set of $q$-terms and  $s$-terms. Then  $(\mathbb{T}_{q,s};\trianglelefteq)$ is bqo, hence $(\mathcal{T}(\mathbb{T}_{q,s});\leq_h)$ is also bqo. The next definition makes precise the relation between the introduced qo's $\trianglelefteq$ and $\leq_h$. 

\begin{Definition}\label{relat}
We associate with any  $u\in\mathbb{T}_\sigma$ the labeled tree $T(u)$ by induction as follows: $T(q)$ is the singleton tree  labeled by $q$, $T(s_\alpha(u))$ is the singleton tree  labeled by $s_\alpha(u)$, $T(F_q(u_0,\ldots))=q\to(T(u_0)\sqcup T(u_1)\sqcup\cdots)$, $T(F_\alpha(u_0,\ldots))=s_\alpha(u_0)\to(T(u_1)\sqcup T(u_2)\sqcup\cdots)$. 
\end{Definition}

Obviously, $T(u)$ is a singleton tree iff $u\in\mathbb{T}_{q,s}$. The next lemma is checked by cases from Definition \ref{treeterms} using induction on terms.

\begin{Lemma}\label{hqo}
The function $u\mapsto T(u)$ is an isomorphism between $(\mathbb{T}_\sigma;\trianglelefteq)$ and $(\mathcal{T}(\mathbb{T}_{q,s});\leq_h)$. 
\end{Lemma}

The next lemma is immediate by induction on terms.

\begin{Lemma}\label{syntax}
Any term $u\in\mathbb{T}_\sigma$ satisfies precisely one of the following alternatives:
\begin{enumerate}
\item $u=q$ for a unique $q\in Q$;
\item $u=s_{\beta_0}\cdots s_{\beta_m}(q)$ for unique $m<\omega$, $\beta_0,\ldots,\beta_m<\omega_1$, $q\in Q$; 
\item $u=F_q(u_0,\ldots)$ for unique $q\in Q$ and $u_0,\ldots\in\mathbb{T}_\sigma$;
\item $u=F_\alpha(u_0,\ldots)$ for unique $\alpha<\omega_1$, and $u_0,\ldots\in\mathbb{T}_\sigma$;
\item $u=s_{\beta_0}\cdots s_{\beta_m}(F_q(u_0,\ldots))$ for unique $m<\omega$, $\beta_0,\ldots,\beta_m<\omega_1$, $q\in Q$,  $u_0,\ldots\in\mathbb{T}_\sigma$;
\item $u=s_{\beta_0}\cdots s_{\beta_m}(F_\alpha(u_0,\ldots))$ for unique $m<\omega$, $\beta_0,\ldots,\beta_m<\omega_1$, $\alpha<\omega_1$,  $u_0,\ldots\in\mathbb{T}_\sigma$.
 \end{enumerate}
 \end{Lemma}
 
Terms from items (1,2) above will be called {\em singleton terms}. With any singleton term $u$ a unique element $q\in Q$ is associated denoted by $q(u)$. Below we will also need the following technical notions.

\begin{Definition}\label{syntax2}
We associate with any $u\in\mathbb{T}_\sigma$ the ordinal $sh(u)$ and the term $u'\in\mathbb{T}_\sigma$ as follows: if $u$ is not an $s$-term then $sh(u)=0$ and $u'=u$, otherwise $sh(u)=\omega^{\beta_0}+\cdots+\omega^{\beta_m}$ and $u'=q$ if $u$ satisfies (2),  $u'=F_q(u_0,\ldots)$ if $u$ satisfies  (5), and $u'=F_\alpha(u_0,\ldots)$ if $u$ satisfies (6).  
 \end{Definition}
 
 We collect some obvious properties of $u'$.

\begin{Lemma}\label{syntax4}
\begin{enumerate}
\item $u'=u$ iff $u$ is not an $s$-term.
\item $u'$ is a subterm of $u$, so $u'\trianglelefteq u$ and if $u$ is an $s$-term then $rk(u')<rk(u)$.
\item $u'$ is not an $s$-term, hence $u''=u$.
\item $u'\in Q$ iff $u$ is a singleton term.
%\item If $u$  is an $s$-term then $rk(u')< rk(u)$.
 \end{enumerate}
 \end{Lemma}
 
\begin{Definition}\label{syntax3}
We associate with any non-singleton term $u\in\mathbb{T}_\sigma$ the set $\mathcal{F}(u)$ of  sequences $S=(\tau_0,\ldots)$ in $\omega^*$ as follows:   
$\tau_0\in T(u')=(T_0,t_0)$;  if $t_0(\tau_0)$  is a singleton term then $S=(\tau_0)$, otherwise $\tau_1\in T(t_0(\tau_0)')=(T_1,t_1)$; if $t_1(\tau_1)$  is a singleton term then $S=(\tau_0,\tau_1)$, otherwise $\tau_2\in T(t_1(\tau_1)')=(T_2,t_2)$; and so on.
 \end{Definition}
 
\begin{Lemma}\label{syntax5}
\begin{enumerate}
\item For any $u\in\mathbb{T}_\sigma$ and $\tau\in T(u)$, $t_u(\tau)\trianglelefteq u$, where $t_u$ is the labeling function on $T(u)$, and $rk(t_u(\tau))\leq rk(u)$.
\item If $u$ is not a singleton term then $rk(t_u(\tau)')< rk(u)$ for every $\tau\in T(u)$.
\item For any non-singleton term $u\in\mathbb{T}_\sigma$, every sequence in $\mathcal{F}(u)$ is finite.
 \end{enumerate}
 \end{Lemma}

{\em Proof.} (1) For $u\in\mathbb{T}_{q,s}$ the assertion is obvious because $\tau=\varepsilon$ and $t_u(\tau)= u$. Let $u=F_q(u_0,\ldots)$, then either $\tau=\varepsilon$ or $\tau\in T(u_i)$ for a unique $i\geq0$. In the first case $t_u(\tau)= q\trianglelefteq u$ and $rk(t_u(\tau))=0< rk(u)$. In the second case by induction we have $t_u(\tau)= t_{u_i}(\tau)\trianglelefteq u_i\trianglelefteq u$ and $rk(t_u(\tau))=rk(t_{u_i}(\tau))\leq rk(u_i)<rk(u)$.

Finally, let $u=F_\alpha(u_0,\ldots)$. Then either $\tau=\varepsilon$ or $\tau\in T(u_i)$ for a unique $i\geq1$. In the first case $t_u(\tau)= s_\alpha(u_0)\trianglelefteq u$ and $rk(t_u(\tau))=rk(s_\alpha(u_0))=rk (u_0)+1\leq rk(u)$. In the second case by induction we have $t_u(\tau)= t_{u_i}(\tau)\trianglelefteq u_i\trianglelefteq u$ and $rk(t_u(\tau))=rk(t_{u_i}(\tau))\leq rk(u_i)<rk(u)$.

(2) Since $u$ is not  singleton, $u$ is not a $q$-term. If $u$ is an $s$-term then $t_u(\tau)=u$, so by Lemma \ref{syntax4}(2) we have $rk(t_u(\tau)')=rk(u')<rk(u)$. If $u=F_q(u_0,\ldots)$ then, by the proof of item (1), $rk(t_u(\tau)')\leq rk(t_u(\tau))<rk(u)$. Finally, let $u=F_\alpha(u_0,\ldots)$. For $\tau\not=\varepsilon$ the assertion follows again from the proof of item (1). For $\tau=\varepsilon$ we have $t_u(\varepsilon)=s_\alpha(u_0)$, hence, by the proof of item (1) and Lemma \ref{syntax4}(2), $rk(t_u(\varepsilon)')=rk(s_\alpha(u_0)')<rk(s_\alpha(u_0))\leq rk(u)$.

(3) Suppose the contrary: the sequence $\tau_0,\tau_1,\ldots$ from Definition \ref{syntax2} is infinite, hence all terms $t_0(\tau_0),t_1(\tau_1),\dots$ are not singleton. By item (2) we then have $rk(u'))> rk(t_0(\tau_0)')> rk(t_1(\tau_1)')>\cdots$, contradicting the well-foundedness of syntactic trees.
 \qed

With any $(\tau_0,\ldots,\tau_m)\in\mathcal{F}(u)$ we associate the constant $q(t_m(\tau_m))\in Q$. For any $q\in Q$ we set $\mathcal{F}_q(u)=\{(\tau_0,\ldots,\tau_m)\in\mathcal{F}(u)\mid q=q(t_m(\tau_m))\}$.

To be more consistent with notation of previous papers and of Introduction,  we sometimes denote $\mathcal{T}(\mathbb{T}_{q,s})$ by $\mathcal{T}_{\omega_1}(Q)$ and use the structures from Lemma \ref{hqo} interchangeably. Let  $\mathcal{T}^\sqcup_{\omega_1}(Q)$ be the set of non-empty labeled forests obtained from trees in $\mathcal{T}_{\omega_1}(Q)$ by deleting the root (alternatively and equivalently, one can think of $\mathcal{T}^\sqcup_{\omega_1}(Q)$ as the set of countable disjoint unions of trees in $\mathcal{T}_{\omega_1}(Q)$). The relation $\leq_h$ is extended to the larger structure of forests in the obvious way.

As observed in \cite{s19}, any qo $Q$ induces a kind of free $\sigma$-semilattice $Q^\sqcup$  which we define as the qo $(Q^*;\leq^*)$ where $Q^*$ is the set of non-empty countable subsets of $Q$ with the so called {\em domination qo} defined by $S\leq^*R$ iff $\forall s\in S\exists r\in R(s\leq_Q r)$. The $\sigma$-join-irreducible elements of $Q^\sqcup$ form an isomorphic copy of $Q$. The operation $\bigsqcup$ of countable supremum in $Q^\sqcup$ is induced by the operation of countable union in $Q^*$. The construction of $\mathcal{T}_{\omega_1}^\sqcup(Q)$ from $\mathcal{T}_{\omega_1}(Q)$ above is a particular case of this general construction. Categorical properties of the construction $Q\mapsto Q^\sqcup$ and characterisations of some algebras expanding $(\mathcal{T}_{\omega^\alpha}^\sqcup(Q);\leq_h)$ as free structures are considered in \cite{s18}.

The characterisation of  $\mathcal{W}_Q$ (see Introduction) in terms of the iterated labeled trees may be now described  as follows (see \cite{km17} for more details). The relation $\simeq$ below denotes the equivalence of qo's.

\begin{Proposition}\cite{km17}\label{char1}
For any countable bqo $Q$, $(\mathcal{T}^\sqcup_{\omega_1}(Q);\leq_h)\simeq(\mathbf{\Delta}^1_1(Q^\mathcal{N});\leq_W)$. The isomorphism of quotient-posets is induced by a map $\mu:\mathcal{T}_{\omega_1}(Q)\to \mathbf{\Delta}^1_1(Q^\mathcal{N})$  sending trees onto the $\sigma$-join irreducible elements.
 \end{Proposition}
 
For more details on the map $\mu$ see Section \ref{qwadgeb} below. There are similar descriptions of natural initial segments of $\mathcal{W}_Q$. 
For any $\gamma<\omega_1$, apply the construction above to the smaller signature $\sigma(Q,\gamma)=\{q,s_\alpha,F_q,F_\alpha\mid q\in Q,\alpha<\gamma\}$ in place of $\sigma(Q,\omega_1)$. The resulting set of labeled trees is denoted by $\mathcal{T}_{\omega^\gamma}(Q)$. We obtain an operator $\mathcal{T}_{\omega^\gamma}$ on BQO.

 Finally, for any $\alpha<\omega_1$ we define the operator $\mathcal{T}_{\alpha}$ on BQO as follows: $\mathcal{T}_0$ is the identity operator, and for any positive countable ordinal $\alpha$ we set $\mathcal{T}_\alpha=\mathcal{T}_{\omega^{\alpha_0}}\circ\cdots\circ\mathcal{T}_{\omega^{\alpha_n}}$ where $n<\omega$ and $\alpha_0\geq\cdots\geq\alpha_n$ are the unique ordinals with $\alpha=\omega^{\alpha_0}+\cdots+\omega^{\alpha_n}$. The set of forests $\mathcal{T}^\sqcup_{\alpha}(Q)$ is obtained from $\mathcal{T}_{\alpha}(Q)$ by the above construction. In particular, $\mathcal{T}_{\alpha+1}=\mathcal{T}_{\alpha}\circ\mathcal{T}$ where $\mathcal{T}$ is the operator from the beginning of this subsection.
 
\begin{Proposition}\cite{km17}\label{char2}
For any countable bqo $Q$ and any $\alpha<\omega_1$ we have: $(\mathcal{T}^\sqcup_{\alpha}(Q);\leq_h)\simeq(\mathbf{\Delta}^0_{1+\alpha}(Q^\mathcal{N});\leq_W)$. 
 \end{Proposition}

We conclude this subsection by a lemma on automorphisms of $(\mathbb{T}_\sigma;\trianglelefteq)$.

\begin{Lemma}\label{auto}
For every poset $(Q;\leq_Q)$, the automorphism group $Aut(Q)$ of $(Q;\leq_Q)$ is isomorphic to a subgroup of the automorphism group $Aut(\mathbb{T}_\sigma)$ of $(\mathbb{T}_\sigma;\trianglelefteq)$. 
 \end{Lemma}
 
 {\em Proof.} We extend any $g\in Aut(Q)$ to a function on $\mathbb{T}_\sigma$ (also denoted by $g$) by induction as follows: $g(q)=g(q)$ for $q\in Q$, $g(s_\alpha(u))=s_\alpha(g(u))$, $g(F_q(u_0),\ldots))=F_{g(q)}(g(u_0),\ldots)$, $g(F_\alpha(u_0,\ldots)=F_\alpha(g(u_0),\ldots)$. By induction, considering 16  cases from Definition \ref{treeterms}, it is straightforward to check that $u\trianglelefteq v$ iff $g(u)\trianglelefteq g(v)$, and $g^{-1}g(u)=u$, for all $u,v\in\mathbb{T}_\sigma$ and $g\in Aut(Q)$.
 \qed

\subsection{Hierarchy bases}\label{bases}

We recall (see e.g. \cite{s16}) the technical notion of a (hierarchy) base. Such bases serve as a starting point for constructing the $Q$-IFH \cite{s16}. They have nothing in common with topological bases.

\begin{Definition}\label{base}
By a base  in a set $X$ we mean a sequence ${\mathcal L}(X)=\{{\mathcal L}_\alpha\}_{\alpha<\omega_1}$, ${\mathcal L}_\alpha={\mathcal L}_\alpha(X)\subseteq P(X)$, such that every ${\mathcal L}_\alpha$ is closed under countable union and finite intersection (in particular, $\emptyset, X\in{\mathcal L}_\alpha$), and ${\mathcal L}_\alpha\cup\check{\mathcal L}_\alpha\subseteq{\mathcal L}_\beta\cap\check{\mathcal L}_\beta$ for all $\alpha<\beta<\omega_1$. 
\end{Definition}

A major natural example of a  hierarchy base in a topological space $X$ is the {\em Borel base} $\mathcal{L}(X)=\{\bfSig^0_{1+\alpha}(X)\}_{\alpha<\omega_1}$. Other natural examples are the  hyperprojective hierarchy and many of its refinements. There are of course many ``unnatural''  bases, e.g. the bases $\{\mathbf{B}(X),\mathbf{B}(X),\ldots\}$ and $\{P(X),P(X),\cdots\}$ over which any IFH of sets collapses to the first level.

With any base $\mathcal{L}(X)$ in $X$ we  associate some new bases as follows. For any $\beta<\omega_1$, let $\mathcal{L}^\beta(X)=\{\mathcal{L}_{\beta+\alpha}(X)\}_{\alpha}$; we call this base in $X$ the {\em $\beta$-shift of $\mathcal{L}(X)$}. For any $U\subseteq X$, let $\mathcal{L}(U)=\{\mathcal{L}_{\alpha}(U)\}$ where $\mathcal{L}_{\alpha}(U)=\{U\cap S\mid S\in\mathcal{L}_{\alpha}(X)\}$; we call this base in $U$ the {\em $U$-restriction of $\mathcal{L}(X)$}. 

\begin{Lemma}\label{baseshift}
\begin{enumerate}
\item $(\mathcal{L}^\beta)^\gamma(X)=\mathcal{L}^{\beta+\gamma}(X)$.
\item If  $\beta*<\alpha*$ (see Section \ref{bqo}) then $\mathcal{L}^\beta_\alpha(X)=\mathcal{L}_\alpha(X)$. Therefore, many levels of $\mathcal{L}(X)$ remain unchanged under the $\beta$-shift.
\end{enumerate}
\end{Lemma}

{\em Proof.} (1) Indeed, $(\mathcal{L}^\beta)^\gamma_\alpha=\mathcal{L}^\beta_{\gamma+\alpha}=\mathcal{L}_{\beta+(\gamma+\alpha)}=\mathcal{L}_{(\beta+\gamma)+\alpha}=\mathcal{L}^{\beta+\gamma}_\alpha$.

(2) Since $\beta+\alpha=\alpha$ by the definition of $\beta*$ and $\alpha*$, $\mathcal{L}^\beta_\alpha(X)=\mathcal{L}_{\beta+\alpha}(X)=\mathcal{L}_\alpha(X)$. 
 \qed

By a {\em morphism} $g:\mathcal{L}(X)\to\mathcal{L}(Y)$ of bases we mean a function $g:P(X)\to P(Y)$ such that $g(\emptyset)=\emptyset$, $g(X)=Y$, $g(\bigcup_nU_n)=\bigcup_ng(U_n)$ for every countable sequence $\{U_n\}$ in $P(X)$ (so, in particular, $U\subseteq V$ implies $g(U)\subseteq g(V)$), and $U\in{\mathcal L}_\alpha(X)$ implies $g(U)\in{\mathcal L}_\alpha(Y)$ for each $\alpha<\omega_1$. Obviously, the identity function on $P(X)$ is a morphism of any base in $X$ to itself, and if $g:\mathcal{L}(X)\to\mathcal{L}(Y)$ and $h:\mathcal{L}(Y)\to\mathcal{L}(Z)$ are morphisms of bases then $h\circ g:\mathcal{L}(X)\to\mathcal{L}(Z)$ is also a morphism. 
We illustrate the notion of morphism with the following well known fact. Recall that a function $f:X\to Y$ between spaces is {\em $\bfSig^0_{1+\alpha}$-measurable} iff $f^{-1}(U)\in\bfSig^0_{1+\alpha}(X)$ for any open set $U$ in $Y$.

\begin{Lemma}\label{measur}
Let $f:X\to Y$ be $\bfSig^0_{1+\alpha}$-measurable and let $\mathcal{L}(X),\mathcal{L}(Y)$ be the Borel bases in $X,Y$ resp. Then $f^{-1}:P(Y)\to P(X)$ is a morphism from $\mathcal{L}(Y)$ to $\mathcal{L}^\alpha(X)$. In particular, if $f$ is continuous then $f^{-1}:P(Y)\to P(X)$ is a morphism of $\mathcal{L}(Y)$ to $\mathcal{L}(X)$. 
\end{Lemma} 

The following class of bases will be frequently mentioned in the sequel.

\begin{Definition}\label{redbase}
A base ${\mathcal L}(X)$ is reducible if every its level ${\mathcal L}_\alpha(X)$ has the $\sigma$-reduction property. 
\end{Definition}

The next fact is known (see e.g. \cite{ke95} and \cite{s13}). 

\begin{Lemma}\label{redbase1}
The Borel base in every zero-dimensional cb$_0$-space is reducible. The 1-shift of the Borel base in every  cb$_0$-space is reducible.
\end{Lemma}

We conclude this subsection with introducing  some auxiliary notions used in the sequel.
For any tree $T\subseteq\omega^*$ and a $T$-family $\{U_\tau\}$ of subsets of $X$, we define the $T$-family $\{\tilde{U}_\tau\}$ of subsets of $X$ by $\tilde{U}_\tau=U\tau\setminus\bigcup\{U_{\tau'}\mid\tau\sqsubset\tau'\in T\}$; the sets $\tilde{U}_\tau$ will be called {\em components} of the family $\{U_\tau\}$. The $T$-family $\{U_\tau\}$ is {\em monotone} if $U_\tau\supseteq U_{\tau'}$ for all $\tau\sqsubseteq\tau'\in T$. We associate with any $T$-family $\{U_\tau\}$ the monotone $T$-family $\{U'_\tau\}$ by $U'_\tau=\bigcup_{\tau'\sqsupseteq\tau}U_{\tau'}$.

\begin{Lemma}\label{tfam}
Let $T$ be a well founded tree, $\mathcal{L}(X)$ be  a base, and $\{U_\tau\}$ be a $T$-family of ${\mathcal L}_\alpha$-sets. Then the components are differences of ${\mathcal L}_\alpha$-sets (hence they belong to ${\mathcal L}_{\alpha+1}\cap\check{\mathcal L}_{\alpha+1}$), $\bigcup_\tau U_\tau=\bigcup_\tau\tilde{U}_\tau$, $\tilde{U}_\tau=\widetilde{U'}_\tau$, and $\tilde{U}_\tau\cap\tilde{U}_{\tau'}=\emptyset$ for $\tau\sqsubset\tau'\in T$.
\end{Lemma}

{\em Proof.} We check only the second assertion, the proofs of others being even simpler. Since $\tilde{U}_\tau\subseteq U_\tau$, $\bigcup_\tau U_\tau\supseteq\bigcup_\tau\tilde{U}_\tau$. Conversely, let $x\in\bigcup_\tau  U_\tau$. Then the set $\{\tau\in T\mid x\in U_\tau\}$ is nonempty. Since $(T;\sqsupseteq)$ is well founded, $ x\in U_\tau$ for some maximal element $\tau$ of $(\{\tau\in T\mid x\in U_\tau\};\sqsubseteq)$; but then $x\in\tilde{U}_\tau$.
 \qed

The next lemma is also easy.
 
\begin{Lemma}\label{tfamunion}
Let $T$ be a well founded tree, $\mathcal{L}(X)$ be  a base,  $\{U^i_\tau\}_i$ be a sequence of monotone $T$-families of ${\mathcal L}_\alpha$-sets, and $U_\tau=\bigcup_iU^i_\tau$ for each $\tau\in T$. Then $\{U_\tau\}$ is a monotone $T$-family of ${\mathcal L}_\alpha$-sets and $\tilde{U}_\tau\subseteq\bigcup_i\widetilde{U^i_\tau}$ for each $\tau\in T$.
\end{Lemma}

We call a $T$-family $\{V_\tau\}$ of ${\mathcal L}_\alpha$-sets {\em reduced} if it is monotone and satisfies $V_{\tau i}\cap V_{\tau j}=\emptyset$ for all $\tau i,\tau j\in T$. Obviously, for any reduced $T$-family $\{V_\tau\}$ of ${\mathcal L}_\alpha$-sets the components $\tilde{V}_\tau$ are pairwise disjoint. The next lemma is checked by a top-down (assuming that  trees grow downwards) application of the $\sigma$-reduction property.

\begin{Lemma}\label{reduce}
Let $T$ be an infinitely-branching well founded tree, $\mathcal{L}(X)$ be  a base,  $\{U_\tau\}$ be a monotone $T$-family of ${\mathcal L}_\alpha$-sets, and let ${\mathcal L}_\alpha$ have the $\sigma$-reduction property. Then there is a reduced $T$-family $\{V_\tau\}$ of ${\mathcal L}_\alpha$-sets such that $V_\tau\subseteq U_\tau$, $\bigcup_\tau V_\tau=\bigcup_\tau U_\tau$, $\bigcup_i\{V_{\tau i}\mid \tau i\in T\}=\bigcup_i\{V_\tau\cap U_{\tau i}\mid \tau i\in T\}$, and $\tilde{V}_\tau\subseteq \tilde{U}_\tau$ for each $\tau\in T$.
\end{Lemma}

{\em Proof.} If $T=\{\varepsilon\}$ is singleton, there is nothing to prove. Otherwise, let $\{V_i\}$ be a reduct of $\{U_i\}$ and let $U'_{i\tau}=V_i\cap U_{i\tau}$ for all $i\tau\in T$. Apply this procedure to the trees $T(i)$ and further downwards whenever possible. Since $T$ is well founded, we will finally obtain a desired reduced family which we  call a reduct of $\{U_\tau\}$.
 \qed

\begin{Lemma}\label{reduce3}
For every well founded tree $T$, a base $\mathcal{L}$(X), $\rho\in T$ and $\alpha<\omega_1$, there is a unique reduced $T$-family  $\{U_\tau\}$ of ${\mathcal L}_\alpha$-sets such that $\tilde{U}_\rho=X$ (and then necessarily $\tilde{U}_\tau=\emptyset$ for all $\tau\in T\setminus\{\rho\}$).
\end{Lemma}

{\em Proof.} Obviously, it is enough to set $U_\tau=X$ if $\tau\sqsubseteq\rho$ and $U_\tau=\emptyset$ otherwise.
 \qed

\subsection{Defining $Q$-partitions by iterated families}\label{itfam}

Here we define the notion of a $u$-family ($u\in\mathbb{T}_\sigma$) in a given base $\mathcal{L}(X)$ and explain how such (iterated) families determine $Q$-partitions of $X$. The definition follows the definition of terms in Section \ref{trees}, induction scheme of Definition \ref{relat} and Lemma \ref{hqo}. 

\begin{Definition}\label{famin}
\begin{enumerate}
\item $F$ is a $q$-family in $\mathcal{L}(X)$ iff $F=\{X\}$.
\item $F$ is an $s_\alpha(u)$-family in $\mathcal{L}(X)$ iff $F$ is a $u$-family in $\mathcal{L}^{\omega^\alpha}(X)$.
\item $F$ is an $F_q(u_0,\ldots)$-family in $\mathcal{L}(X)$ iff $F=(\{U_\tau\},\{F_\tau\})$ where $\{U_\tau\}$ is a monotone $T$-family of  ${\mathcal L}_0$-sets with $U_\varepsilon=X$ and, for each $\tau\in T$, $F_\tau$ is  a  $t(\tau)$-family in $\mathcal{L}(\tilde{U}_\tau)$, where $(T,t)=T(F_q(u_0,\ldots))$.
\item $F$ is an $F_\alpha(u_0,\ldots)$-family in $\mathcal{L}(X)$ iff $F=(\{U_\tau\},\{F_\tau\})$ where $\{U_\tau\}$ is a monotone $T$-family of  ${\mathcal L}_0$-sets with $U_\varepsilon=X$ and, for each $\tau\in T$, $F_\tau$ is  a  $t(\tau)$-family in $\mathcal{L}(\tilde{U}_\tau)$, where $(T,t)=T(F_\alpha(u_0,\ldots))$.
\end{enumerate}
\end{Definition}

The notion of a reduced $u$-family $F$ is obtained from this definition by requiring $\{U_\tau\}$ and $F_\tau$ in items (3,4) to be reduced. Note that Definition \ref{famin} and the next definition are uniform in bases, i.e., for any fixed $u$, the  $u$-family $F$ in any item above is defined for all bases simultaneously.

From Lemma \ref{syntax} we obtain the following information on the structure of $u$-families in $\mathcal{L}(X)$ where we use notions from Definition \ref{syntax2}.

\begin{Lemma}\label{syntax1}
Let $F$ be a $u$-family in $\mathcal{L}(X)$.
If $u$ is a singleton term then $F=\{X\}$, otherwise $F=(\{U_\tau\},\{F_\tau\})$ where $\{U_\tau\}$ is a monotone $T(u')$-family of  ${\mathcal L}^{sh(u)}_0$-sets with $U_\varepsilon=X$ and, for each $\tau\in T(u')$, $F_\tau$ is  a  $t(\tau)$-family in $\mathcal{L}^{sh(u)}(\tilde{U}_\tau)$.
 \end{Lemma}

Now we define the notion ``a  $u$-family $F$ in $\mathcal{L}(X)$ determines a partition $A:X\to Q$''. In general, every $u$-family determines at most one $Q$-partition, not every $u$-family determines a $Q$-partition, and every reduced $u$-family determines a $Q$-partition.

\begin{Definition}\label{determ}
\begin{enumerate}
\item A $q$-family $F$ in $\mathcal{L}(X)$ determines $A$ iff $A=\lambda x.q$.
\item An $s_\alpha(u)$-family $F$ in $\mathcal{L}(X)$ determines $A$ iff $F$ determines $A$ as a $u$-family in $\mathcal{L}^{\omega^\alpha}(X)$.
\item For $u\in\{F_q(u_0,\ldots),F_\alpha(u_0,\ldots)\}$, a $u$-family $F=(\{U_\tau\},\{F_\tau\})$ in $\mathcal{L}(X)$  determines $A$ iff
for each $\tau\in T(u)$, $F_\tau$ determines $A|_{\tilde{U}_\tau}$. 
\end{enumerate}
\end{Definition}

By definitions above and Lemma \ref{syntax1}, a $u$-family $F$ in $\mathcal{L}(X)$ that determines a $Q$-partition $A$ yields a mind-change ``algorithm'' for computing $A(x)$ for any given $x\in X$ as follows. We use the set $\mathcal{F}(u)$ from Definition \ref{syntax3} and Lemma \ref{syntax5}.
 
If $u$ is a singleton term, $A$ is the constant $Q$-partition $\lambda x.q(u)$, hence $A(x)=q(u)$. Otherwise, $F=(\{U_{\tau_0}\},\{F_{\tau_0}\})$ where $\{U_{\tau_0}\}$ is a monotone $u'$-family of  ${\mathcal L}^{sh(u)}_0$-sets with $U_\varepsilon=X$ and, for each $\tau_0\in T(u')$, $F_{\tau_0}$ is  a  $t_0(\tau_0)$-family in $\mathcal{L}^{sh(u)}(\tilde{U}_{\tau_0})$ (which coincides with the $t_0(\tau_0)'$-family in $\mathcal{L}^{sh(u)+sh(t_0(\tau_0))}(\tilde{U}_{\tau_0})$).  Since the components  $\tilde{U}_{\tau_0}$ (which we call first level components of $F$) cover $X$ by Lemma \ref{tfam}, $x\in\tilde{U}_{\tau_0}$ for some $\tau_0\in T(u')$; $\tau_0$ is searched by the usual mind-change procedure working with differences of ${\mathcal L}^{sh(u)}_0$-sets (see Lemma \ref{tfam}).

If the term $t_0(\tau_0)$ is singleton, $A|_{\tilde{U}_{\tau_0}}$ is a constant $Q$-partition and we have computed $A(x)\in Q$. Otherwise, $F_{\tau_0}=(\{U_{\tau_0\tau_1}\},\{F_{\tau_0\tau_1}\})$ and we can continue the computation as above and find a second level component $\tilde{U}_{\tau_0\tau_1}$ of $F$ containing $x$; this is a harder  mind-change procedure working with differences of ${\mathcal L}^{sh(u)+sh(t_0(\tau_0))}_0$-sets. We continue this process until we reach a sequence  $(\tau_0,\ldots,\tau_m)\in\mathcal{F}(u)$ such that $x\in\tilde{U}_{\tau_0\cdots\tau_m}$ and $t_m(\tau_m)$ is a singleton term; such components $\tilde{U}_{\tau_0\cdots\tau_m}$ are called {\em terminating} and have the associated constants $q(\tau_0,\ldots,\tau_m)=q(t_m(\tau_m))\in Q$. Note that the terminating components cover $X$ and if the family $F$ is reduced then the terminating components form a partition of $X$. In any case we have: $A^{-1}(q)=\bigcup\{\tilde{U}_{\tau_0\cdots\tau_m}\mid (\tau_0,\ldots,\tau_m)\in\mathcal{F}_q(u)\}$ for each $q\in Q$.

Note also that if the family $F$ above was reduced then the computation is ``linear'' since the components of each level are pairwise disjoint and cover the parent component, otherwise the computation is ``parallel'' since already at the first level $x$ may belong to several components $\tilde{U}_{\tau_0}$. 

The described procedure enables to write a $u$-family $F$, where $u$ is not a singleton term, in an explicit (but not completely precise) form of $u'$-family $(\{U_{\tau_0}\},\{U_{\tau_0\tau_1}\},\ldots)$ in  ${\mathcal L}^{sh(u)}(X)$  which is sometimes more intuitive than the form $(\{U_{\tau}\},\{F_\tau\})$ above.

We formulate some properties of the introduced notions. The next lemma is immediate by definitions. 

\begin{Lemma}\label{ui}
Let $u$ be a non-singleton term and the $u'$-family $(\{U_{\tau_0}\},\{U_{\tau_0\tau_1}\},\ldots)$ in  ${\mathcal L}^{sh(u)}(X)$ determines $A\in Q^X$. 
\begin{enumerate}
\item If $u'=F_q(u_0,\ldots)$ then $A|_{U_i}$ is determined by the $u_i$-family $(\{U_{i\sigma_0}\},\{U_{i\sigma_0\tau_1}\},\ldots)$ in  ${\mathcal L}^{sh(u)}(U_i)$, for each $i\geq0$.
\item If $u'=F_\alpha(u_0,\ldots)$ then $A|_{U_{i+1}}$ is determined by the $u_{i+1}$-family $(\{U_{i\sigma_0}\},\{U_{i\sigma_0\tau_1}\},\ldots)$ in  ${\mathcal L}^{sh(u)}(U_{i+1})$, for each $i\geq0$.
\end{enumerate}
 \end{Lemma}

Let $f:X\to Y$ be a function such that $f^{-1}$ is a morphism from $\mathcal{L}(Y)$ to $\mathcal{L}(X)$. Associate with any $u$-family $F$ in $\mathcal{L}(Y)$ the $u$-family $f^{-1}(F)$ in $\mathcal{L}(X)$ as follows:  if $u=q$ then $f^{-1}(F)=\{X\}$; if $u=s_\alpha(v)$ then $f^{-1}(F)$ is the $v$-family $f^{-1}(F)$ in $\mathcal{L}^{\omega^\alpha}(X)$; in the remaining cases we have $F=(\{U_\tau\},\{F_\tau\})$, and we set $f^{-1}(F)=(\{f^{-1}(U_\tau)\},\{f^{-1}(F_\tau)\})$. Obviously, $f^{-1}(F)$ is indeed a $u$-family in $\mathcal{L}(X)$. The next lemma is immediate by induction.

\begin{Lemma}\label{preim}
In assumptions of the previous paragraph, if a $u$-family $F$ in $\mathcal{L}(Y)$ determines $A$ then the $u$-family $f^{-1}(F)$ in $\mathcal{L}(X)$ determines $A\circ f$.
 \end{Lemma}

Now we associate with any $u$-family $F$ in $\mathcal{L}(X)$ and any $V\subseteq X$ the $u$-family $F|_V$ in the $V$-restriction $\mathcal{L}(V)$ (see Section \ref{bases}) as follows:  if $u=q$ then $F|_V=\{V\}$; if $u=s_\alpha(v)$ then $F|_V$ is the $v$-family $F|_V$  in $\mathcal{L}^{\omega^\alpha}(V)$; in the remaining cases we have $F=(\{U_\tau\},\{F_\tau\})$, and we set $F|_V=(\{V\cap U_\tau\},\{F_\tau|_V\})$. Obviously, $F|_V$ is indeed a  $u$-family in $\mathcal{L}(V)$. The next lemma is immediate by induction.

\begin{Lemma}\label{restrict}
In assumptions of the previous paragraph, if a $u$-family $F$ in $\mathcal{L}(X)$  determines $A$ then the $u$-family $F|_V$ in $\mathcal{L}(V)$ determines $A|_V$.
 \end{Lemma}

Let $\{G_i\}$, $G_i=(\{U^i_{\tau_0}\},\{U^i_{\tau_0\tau_1}\},\ldots)$, be a sequence of $u$-families ($u$ is a non-singleton term) in  $\mathcal{L}(Y_i)$, $Y_i\subseteq X$, then $G=(\{U_{\tau_0}\},\{U_{\tau_0\tau_1}\},\ldots)$, where  $U_{\tau_0}=\bigcup_iU^i_{\tau_0}$, $U_{\tau_0\tau_1}=\bigcup_iU^i_{\tau_0\tau_1}\ldots$, is a $u$-family  in  $\mathcal{L}(Y)$ where $Y=\bigcup_iY_i$.  The next lemma follows from definitions and Lemma \ref{tfamunion}.
 
 \begin{Lemma}\label{famunion}
 Let $A\in Q^X$. 
In assumptions of the previous paragraph, if the $u$-family $G_i$ in $\mathcal{L}(Y_i)$  determines $A|_{Y_i}$ for each $i\geq0$ then the $u$-family $G$ in $\mathcal{L}(Y)$ determines $A|_Y$.
  \end{Lemma}

The next lemma is also clear.
  
   \begin{Lemma}\label{famconstruct}
 Let $A\in Q^X$, $Y\in\mathcal{L}_0(X)\cap\check{\mathcal{L}}_0(X)$, $A(x)=q$ for $x\in X\setminus Y$, let $A|_Y$ be determined by a $u$-family $F$ in $\mathcal{L}(Y)$, and let $\tilde{U}_{\tau_0\cdots\tau_m}$ be a terminating component of $F$ with   $q=q(\tau_0,\ldots,\tau_m)$. Then there is a $u$-family $F'$ in $\mathcal{L}(X)$ such that its $(\tau_0,\ldots,\tau_m)$-terminating component is  $\tilde{U}_{\tau_0\cdots\tau_m}\cup\bar{Y}$, all other terminating components coincide with those of $F$, and $F'$ determines   $A$.
  \end{Lemma}

Let $F=(\{U_{\tau_0}\},\{U_{\tau_0\tau_1}\},\ldots)$ and $G=(\{V_{\tau_0}\},\{V_{\tau_0\tau_1}\},\ldots)$ be  $u$-families in $\mathcal{L}(X)$. We say that $G$ {\em is a reduct of $F$} if $G$ is  reduced and $\tilde{V}_{\tau_0\cdots\tau_m}\subseteq\tilde{U}_{\tau_0\cdots\tau_m}$ for each $(\tau_0,\ldots,\tau_m)\in\mathcal{F}(u)$. 

\begin{Lemma}\label{reduct}
Let ${\mathcal L}(X)$ be a reducible base in $X$ and $u\in\mathbb{T}_\sigma$. Then any $u$-family $F$ in ${\mathcal L}(X)$ has a reduct $G$. Moreover, if $F$ determines $A$ then any reduct of $F$ also determines $A$.
\end{Lemma}

{\em Proof Sketch.} We follow the procedure of computing $A(x)$ described above. If  $u$ is a singleton term, we set $G=F=\{X\}$; then $F,G$ determine the same constant $Q$-partition.  Otherwise, $F$ has the form as above. Let  $G$ as above be obtained from $F$ by repeated reductions from Lemma \ref{reduce}, so in particular $\tilde{V}_{\tau_0\cdots\tau_m}\subseteq\tilde{U}_{\tau_0\cdots\tau_m}$ for each $(\tau_0,\ldots,\tau_m)\in\mathcal{F}(u)$.

For the second assertion, let $F$ determine $A$ and let $G$ be a reduct of $F$. For any $x\in X$, let $\tilde{V}_{\tau_0\cdots\tau_m}$ be the unique terminating component of $G$ containing $x$. Then also $x\in\tilde{U}_{\tau_0\cdots\tau_m}$, hence $A(x)=q(\tau_0,\ldots,\tau_m)$ and $G$ determines $A$.
 \qed
 
The next lemma  follows from the results above.
 
\begin{Lemma}\label{tfaml}
 Every $u$-family $F$ in ${\mathcal L}(X)$ determines at most one $Q$-partition of $X$.
 Every reduced $u$-family $G$ in ${\mathcal L}(X)$ determines precisely one $Q$-partition of $X$.
\end{Lemma}

{\em Proof.} The second assertion follows from the remark that the terminating components of $G$ form a partition of $X$. For the first assertion, let  $F$ in ${\mathcal L}(X)$ determine $Q$-partitions $A,B$ of $X$. Let $x\in X$. If  $u$  is a singleton term, $F$ determines a constant $Q$-partition, so in particular $A(x)=B(x)$.  Otherwise, $F=(\{U_\tau\},\{F_\tau\})$ as specified above. By the procedure of computing $A(x)$, there is a terminating component $\tilde{U}_{\tau_0\cdots\tau_m}\ni x$ of $F$. By Definition \ref{determ}, $A(x)=q(\tau_0,\ldots,\tau_m)=B(x)$.
 \qed

\subsection{Infinitary fine hierarchy over a base}\label{qfineb}

Here we define the  $Q$-IFH over a given base and prove some of its  properties. 

\begin{Definition}\label{fhoverabase}
Associate with any  base $\mathcal{L}(X)$  in $X$, any qo $Q$, and any $u\in\mathbb{T}_\sigma$ the set $\mathcal{L}(X,u)$ of $Q$-partitions of $X$ determined by some $u$-family in $\mathcal{L}(X)$.  The family $\{\mathcal{L}(X,u)\}_{u\in\mathbb{T}_\sigma}$ is called the {\em infinitary $Q$-fine hierarchy over $\mathcal{L}(X)$}. 
\end{Definition}

The algorithm of computing $A(x)$, where $A\in\mathcal{L}(X,u)$ is determined by a $u$-family, described in the preceding subsection, explains in which sence the $Q$-IFH over $\mathcal{L}(X)$ may be considered as an ``iterated difference hierarchy''.

By Lemma \ref{hqo}, we can equivalently denote the $Q$-IFH over $\mathcal{L}(X)$ as $\{\mathcal{L}(X,T)\}_{T\in\mathcal{T}_{\omega_1}(Q)}$, as we did in the Introduction; so now we have precise definitions of the objects discussed there.
The next property describes the behaviour of $Q$-IFH w.r.t. the operations on  bases from Section \ref{bases}.

\begin{Lemma}\label{shift}
\begin{enumerate}
\item For any $\alpha<\omega_1$, $\mathcal{L}(X,s_\alpha(u))=\mathcal{L}^{\omega^\alpha}(X,u)$ and  $\mathcal{L}(X,u)=\mathcal{L}^{sh(u)}(X,u')$. 
 \item For any $V\subseteq X$,  $A\in\mathcal{L}(X,u)$ implies $A|_V\in\mathcal{L}(V,u)$.
\item Let $u$ be  non-singleton  and let $A$  be determined by a $u$-family $(\{U_{\tau_0}\},\{U_{\tau_0\tau_1}\},\ldots)$ in $\mathcal{L}(X)$. If $u'=F_q(u_0,\ldots)$ (resp.  $u'=F_\alpha(u_0,\ldots)$) then $A|_{U_i}\in\mathcal{L}(X,u_i)$ for each $i\geq0$ (resp. $i\geq1$).
 \item Let $A\in Q^X$, $u_0,u_1,\ldots\in\mathbb{T}_\sigma$, and let $\{U_i\}_{i\geq0}$ be non-empty open sets  not exhausting $X$ such that $A|_V=\lambda v.q$ (where $V=\calN\setminus\bigcup_iU_i$) and $A|_{U_i}\in\mathcal{L}(U_i,u_i)$ for all $i\geq0$. Then $A\in\mathcal{L}(X,u)$ where $u=F_q(u_0,\ldots)$. 
 \item Let $A\in Q^X$, $u_0,u_1,\ldots\in\mathbb{T}_\sigma$, and let $\{U_i\}_{i\geq1}$ be non-empty open sets  not exhausting $X$ such that $A|_V\in\mathcal{L}(X,s_\alpha(u_0))$ (where $V=\calN\setminus\bigcup_{i\geq1}U_i$) and $A|_{U_i}\in\mathcal{L}(U_i,u_i)$ for all $i\geq1$. Then $A\in\mathcal{L}(X,u)$ where $u=F_\alpha(u_0,\ldots)$. 
  \end{enumerate}
\end{Lemma}

{\em Proof.} (1) The second assertion follows from the first one which holds by Definition \ref{determ}.

(2) Let $A\in\mathcal{L}(X,u)$ be determined by a $u$-family $F$ in $\mathcal{L}(X)$. By Lemma \ref{restrict}, $A|_V$ is determined by the $u$-family $F|_V$ in $\mathcal{L}(X)$, hence $A|_V\in\mathcal{L}(V,u)$. 

(3) Follows from Lemma \ref{ui}.

(4) Let $A|_{U_i}\in\mathcal{L}(X,u_i)$ be determined by a $u_i$-family $G_i=(\{U^i_{\tau_0}\},\{U^i_{\tau_0\tau_1}\},\ldots)$ in $\mathcal{L}(U_i)$, for each $i\geq0$. By Definition \ref{relat}, $T(u)=q\to(T(u_0)\sqcup T(u_1)\sqcup\cdots)$. We define the $u$-family $G=(\{V_{\tau_0}\},\{V_{\tau_0\tau_1}\},\ldots)$  in $\mathcal{L}(X)$ as follows: $V_\varepsilon=X$, $V_{i\tau_0}=U^i_{\tau_0}$, $V_{i\tau_0\tau_1}=U^i_{\tau_0\tau_1}$, and so on. Then $G$ determines $A$, hence $A\in\mathcal{L}(X,u)$.

(5) Similar to (4).
 \qed

Now we discuss inclusions of levels of the $Q$-IFH.

\begin{Lemma}\label{inclus}
\begin{enumerate}
\item $\mathcal{L}(X,u)\subseteq\mathcal{L}(X,s_\alpha(u))$. 
 \item $\mathcal{L}(X,q)\subseteq\mathcal{L}(X,F_q(u_0,\ldots))$.
  \item $\mathcal{L}(X,u_i)\subseteq\mathcal{L}(X,F_q(u_0,\ldots))$ for all $i\geq0$.
 \item $\mathcal{L}(X,s_\alpha(u_0))\subseteq\mathcal{L}(X,F_\alpha(u_0,\ldots))$.
 \item $\mathcal{L}(X,u_{i+1})\subseteq\mathcal{L}(X,F_\alpha(u_0,\ldots))$ for all $i\geq0$.
 \item Let $u,v\in\mathbb{T}_\sigma$, $\beta,\gamma<\omega_1$, and $\mathcal{L}^\beta(X,u)\subseteq\mathcal{L}^\gamma(X,v)$ over all bases $\mathcal{L}(X)$ in $X$. Then $\mathcal{L}^{\alpha+\beta}(X,u)\subseteq\mathcal{L}^{\alpha+\gamma}(X,v)$ for any $\alpha<\omega_1$ and any base $\mathcal{L}(X)$ in $X$.
\end{enumerate}
\end{Lemma}

{\em Proof.} (1) Let $A\in\mathcal{L}(X,u)$, then $A$ is determined by a $u$-family $F$ in $\mathcal{L}(X)$. By Definition \ref{base}, $F$ is also a $u$-family in $\mathcal{L}^{\omega^\alpha}(X)$, hence $A\in\mathcal{L}^{\omega^\alpha}(X,u)$. By Lemma \ref{shift}(1), $A\in\mathcal{L}(X,s_\alpha(u))$.

(2) We have to show that $\lambda x.q$ is determined by a $u$-family $F=(\{U_\tau\},\{F_\tau\})$ in $\mathcal{L}(X)$, where $u=F_q(u_0,\ldots)$ and $\tau\in T(u)$. Let $\{U_\tau\}$ be the reduced family of $\mathcal{L}_0$-sets with $\tilde{U}_\varepsilon=X$ from Lemma \ref{reduce3}. Let $F_\varepsilon=\{X\}$. For any $\tau\in T(u)\setminus\{\varepsilon\}$, let $F_\tau$ be the trivial reduced $t(\tau)$-family in $\mathcal{L}(\emptyset)$ with empty components. By Definition \ref{relat}, the family $F$ determines $\lambda x.q$. 

(3) Let $A$ be determined by a $u_i$-family $G$ in $\mathcal{L}(X)$. We have to show that $A$ is determined by a $u$-family $F=(\{U_\tau\},\{F_\tau\})$ in $\mathcal{L}(X)$, where $u=F_q(u_0,\ldots)$ and $\tau\in T(u)$. Let $\{U_\tau\}$ be the reduced family of $\mathcal{L}_0$-sets with $\tilde{U}_i=X$ from Lemma \ref{reduce3}. Let $F_i=G$. For any $\tau\in T(u)\setminus\{i\}$, let $F_\tau$ be the trivial reduced $t(\tau)$-family in $\mathcal{L}(\emptyset)$  with empty components. By Definition \ref{relat}, the family $F$ determines $A$.

Items (4,5) are checked by  manipulations similar to those in (2,3). 

(6) For the base $\mathcal{L}^\alpha(X)$ in $X$ the given inclusion reads $(\mathcal{L}^\alpha)^\beta(X,u)\subseteq(\mathcal{L}^\alpha)^\gamma(X,v)$. By Lemma \ref{baseshift}(1), $\mathcal{L}^{\alpha+\beta}(X,u)\subseteq\mathcal{L}^{\alpha+\gamma}(X,v)$.
 \qed

The main result about inclusions of levels of the $Q$-IFH is the following assertion checked by  induction on the 16 cases of Definition \ref{treeterms}, using lemmas above.

\begin{Theorem}\label{inclus1}
If $Q$ is antichain  and $u\trianglelefteq v$, then $\mathcal{L}(X,u)\subseteq\mathcal{L}(X,v)$ for all bases $\mathcal{L}(X)$.
 \end{Theorem}

{\em Proof.} (1) Let $q\trianglelefteq r$, then $q\leq_Qr$, hence $q=r$, hence trivially $\mathcal{L}(X,q)\subseteq\mathcal{L}(X,r)$. 

(2) Let $q\trianglelefteq s_\alpha(u)$, then $q\trianglelefteq u$, hence by induction and Lemma \ref{inclus}(1) $\mathcal{L}(X,q)\subseteq\mathcal{L}(X,u)\subseteq\mathcal{L}(X,s_\alpha(u))$.

(3) Let $q\trianglelefteq F_r(u_0,\ldots)$, then $q\trianglelefteq r$ or $q\trianglelefteq u_i$ for some $i\geq0$, and the inclusion follows by induction and Lemma \ref{inclus}(2,3).

(4) Let $q\trianglelefteq F_\alpha(u_0,\ldots)$, then $q\trianglelefteq s_\alpha(u_0)$ or $q\trianglelefteq u_i$ for some $i\geq1$, and the inclusion follows by induction and Lemma \ref{inclus}(4,5).

(5) Let $s_\alpha(u)\trianglelefteq r$, then $u\trianglelefteq r$. By induction, $\mathcal{L}(X,u)\subseteq\mathcal{L}(X,r)=\{\lambda x.r\}$. By the uniformity of Definition \ref{determ}, $\mathcal{L}(X,s_\alpha(u))=\{\lambda x.r\}$.

(6) Let $s_\alpha(u)\trianglelefteq s_\beta(v)$. Then ($\alpha<\beta$ and $u\trianglelefteq s_\beta(v)$) or ($\alpha=\beta$ and $u\trianglelefteq v$) or ($\alpha>\beta$ and $s_\alpha(u)\trianglelefteq v$). In the first case, by induction we have $\mathcal{L}(X,u)\subseteq\mathcal{L}(X,s_\beta(v))\subseteq\mathcal{L}^{\omega^\beta}(X,v)$. By Lemmas \ref{inclus}(6), \ref{baseshift}(2) and \ref{shift}(1), $\mathcal{L}(X,s_\alpha(u))=\mathcal{L}^{\omega^\alpha}(X,u)\subseteq\mathcal{L}^{\omega^\alpha+\omega^\beta}(X,v)=\mathcal{L}^{\omega^\beta}(X,v)=\mathcal{L}(X,s_\beta(v))$. In the second case, by induction we have $\mathcal{L}(X,u)\subseteq\mathcal{L}(X,v)$, hence $\mathcal{L}^{\omega^\alpha}(X,u)\subseteq\mathcal{L}^{\omega^\beta}(X,v)$, hence $\mathcal{L}(X,s_\alpha(u))\subseteq\mathcal{L}(X,s_\beta(v))$. The third case is even easier.

(7) Let $s_\alpha(u)\trianglelefteq F_r(v_0,\ldots)$, then $s_\alpha(u)\trianglelefteq r$ or $s_\alpha(u)\trianglelefteq v_i$ for some $i\geq0$. The assertion follows by Lemma \ref{inclus}(2) or (3), resp.

(8) Let $s_\alpha(u)\trianglelefteq F_\beta(v_0,\ldots)$, then $s_\alpha(u)\trianglelefteq s_\beta(v_0)$ or $s_\alpha(u)\trianglelefteq v_i$ for some $i\geq1$. The assertion follows by Lemma \ref{inclus}(4) or (5), resp.

(9) Let $F_q(u_0,\ldots)\trianglelefteq r$, then $q\trianglelefteq r$ and $u_i\trianglelefteq r$ for all $i\geq0$. In this case the argument of item (5) works.

(10) Let $F_q(u_0,\ldots)\trianglelefteq s_\alpha(v)$, then $q\trianglelefteq s_\alpha(v)$ and $u_i\trianglelefteq s_\alpha(v)$ for all $i\geq0$. If $v$ is a singleton term, the argument of item (9) works, so let $v$ be a non-singleton term. Without loss of generality we way think that $v$ is an $F$-term (otherwise, $\mathcal{L}^{\omega^\alpha}(X,v)=\mathcal{L}^{\omega^\alpha+sh(v)}(X,v')$, and we can work with the $F$-term $v'$ instead of $v$). 

Let $A\in\mathcal{L}(X,F_q(u_0,\ldots))$, we have to show that $A\in\mathcal{L}(X,s_\alpha(v))$. Let $(\{U_{\tau_0}\},\{U_{\tau_0\tau_1}\},\ldots)$ be a $u$-family  in  $\mathcal{L}(X)$ that determines $A$, then $A(x)=q$ for each $x\in\tilde{U}_\varepsilon$ 
(note that $\tilde{U}_\varepsilon\in\mathcal{L}^{\omega^\alpha}_0(X)\cap\check{\mathcal{L}}^{\omega^\alpha}_0(X)$) and, by Lemma \ref{ui}, $A|_{U_i}$ is determined by the $u_i$-family $(\{U_{i\tau_1}\},\ldots)$ in  $\mathcal{L}(U_i)$ for every $i\geq0$. 
By induction, $A|_{U_i}\in\mathcal{L}^{\omega^\alpha}(U_i,v)$ for every $i\geq0$, so let $G_i=(\{V^i_{\tau_0}\},\{V^i_{\tau_0\tau_1}\},\ldots)$ be a $v$-family  in  $\mathcal{L}^{\omega^\alpha}(U_i)$ that determines $A|_{U_i}$. 
By Lemma \ref{famunion}, the $v$-family $G=\bigcup_iG_i=(\{V_{\tau_0}\},\{V_{\tau_0\tau_1}\},\ldots)$ in $\mathcal{L}^{\omega^\alpha}(\bigcup_iU_i)$   determines    $A_{\bigcup_iU_i}$.
By Lemma \ref{famconstruct}, the $s_\alpha(v)$-family $G'$ determines $A$, hence $A\in\mathcal{L}(X,s_\alpha(v))$.

(11) Let $F_q(u_0,\ldots)\trianglelefteq F_r(v_0,\ldots)$, then ($q\trianglelefteq r$ and $u_i\trianglelefteq F_r(v_0,\ldots)$ for all $i\geq0$) or $F_q(u_0,\ldots)\trianglelefteq v_i$ for some $i\geq0$; the second case follows from Lemma \ref{inclus}(3), so consider the first case. Since $Q$ is antichain, $q=r$.
Let $A\in\mathcal{L}(X,F_q(u_0,\ldots))$, we have to show that $A\in\mathcal{L}(X,F_q(v_0,\ldots))$. Let $(\{U_{\tau_0}\},\{U_{\tau_0\tau_1}\},\ldots)$ be a $u$-family  in  $\mathcal{L}(X)$, where $u=F_q(u_0,\ldots)$, that determines $A$, then $A(x)=q$ for each $x\in\tilde{U}_\varepsilon$, and, by Lemma \ref{ui}, $A|_{U_i}$ is determined by the family $(\{U_{i\tau_1}\},\ldots)$ in  $\mathcal{L}(U_i)$ for each $i\geq0$. 
By induction, $A|_{U_i}\in\mathcal{L}(U_i,v)$ for each $i\geq0$, where  $v=F_q(v_0,\ldots)$, so $A|_{U_i}$ is determined by a $v$-family $G_i=(\{V^i_{\tau_0}\},\{V^i_{\tau_0\tau_1}\},\ldots)$  in  $\mathcal{L}(U_i)$.  
By Lemma \ref{famunion}, the $v$-family $G=(\{V_{\tau_0}\},\{V_{\tau_0\tau_1}\},\ldots)$ in $\mathcal{L}(\bigcup_iU_i)$   determines    $A|_{\bigcup_iU_i}$.
Correcting the $v$-family $G$ by changing $V_\varepsilon$ to $X$, we obtain a $v$-family $G'$ in $\mathcal{L}(X)$ that  determines    $A$. Thus, $A\in\mathcal{L}(X,v)$.  

Items (12,15,16) are checked similar to (10,11), item (13) similar to (9), item (14) similar to (11).
 \qed  

%(12) Let $F_q(u_0,\ldots)\trianglelefteq F_\beta(v_0,\ldots)$ iff ($q\trianglelefteq s_\beta(v_0)$ and $u_i\trianglelefteq F_\beta(v_0,\ldots)$ for all $i\geq0$) or $F_p(u_0,\ldots)\trianglelefteq v_i$ for some $i\geq1$;

%(13) Let $F_\alpha(u_0,\ldots)\trianglelefteq r$ iff $s_\alpha(u_0)\trianglelefteq r$ and $u_i\trianglelefteq q$ for all $i\geq1$;

%(14) Let $F_\alpha(u_0,\ldots)\trianglelefteq s_\beta(v)$ iff $s_\alpha(u_0)\trianglelefteq s_\beta(v)$ and $u_i\trianglelefteq s_\beta(v)$ for all $i\geq1$;

%(15) Let $F_\alpha(u_0,\ldots)\trianglelefteq F_r(v_0,\ldots)$ iff ($s_\alpha(u_0)\trianglelefteq r$ and $u_i\trianglelefteq F_q(v_0,\ldots)$ for all $i\geq1$) or $F_\alpha(u_0,\ldots)\trianglelefteq v_i$ for some $i\geq0$; 

%(16) Let $F_\alpha(u_0,\ldots)\trianglelefteq F_\beta(v_0,\ldots)$ iff ($s_\alpha(u_0)\trianglelefteq s_\beta(v_0)$ and $u_i\trianglelefteq F_\beta(v_0,\ldots)$ for all $i\geq1$) or $F_\alpha(u_0,\ldots)\trianglelefteq v_i$ for some $i\geq1$.

\begin{Corollary}\label{inclus2}
The levels of $\bar{k}$-IFH over any base $\mathcal{L}(X)$ are bqo under inclusion, i.e. for $Q=\bar{k}$ the poset $(\{\mathcal{L}(X,u)\mid u\in\mathbb{T}_\sigma\};\subseteq)$ is bqo. 
 \end{Corollary}

{\em Proof.} By Theorem \ref{inclus1}, $u\mapsto\mathcal{L}(X,u)$ is a monotone surjection from bqo $(\mathbb{T}_\sigma;\trianglelefteq)$ onto $(\{\mathcal{L}(X,u)\mid u\in\mathbb{T}_\sigma\};\subseteq)$. Hence, the latter structure is also bqo.
 \qed

We conclude this subsection with a result about the reduction property. Let the classes
red-$\mathcal{L}(X,u)$ be defined as the classes $\mathcal{L}(X,u)$ in Definition \ref{fhoverabase} but with the reduced families in place of arbitrary families.

\begin{Proposition}\label{red-hier}
If  $\mathcal{L}(X)$ is a reducible base then $\mathcal{L}(X,u)=$red-$\mathcal{L}(X,u)$ for each $u\in\mathbb{T}_\sigma$.
 \end{Proposition}

{\em Proof.} The inclusion from right to left is obvious. Conversely, let $A\in\mathcal{L}(X,u)$, then $A$ is determined by a $u$-family $F$ in $\mathcal{L}(X)$. By Lemma \ref{reduct}, $A$ is determined by a $u$-family $G$ in $\mathcal{L}(X)$ which is a reduct of $F$. Thus, $A$ is  in red-$\mathcal{L}(X,u)$.
 \qed

\section{Infinitary fine hierarchies in cb$_0$-spaces}\label{char}

In this section we study the $Q$-IFH in cb$_0$-spaces. We show that some important properties are preserved by  continuous open surjections while others are not, and we give the set-theoretic description of the $Q$-Wadge hierarchy in  the Baire space. From now on all bases we discuss are the Borel bases $\mathcal{L}(X)=\{\bfSig^0_{1+\alpha}(X)\}_{\alpha<\omega_1}$ in cb$_0$-spaces $X$.

\subsection{General properties}\label{general}

Here we collect some general properties of $Q$-IFH in cb$_0$-spaces. Let $\mathcal{L}(X),\mathcal{L}(Y)$ be the Borel bases in cb$_0$-spaces  $X,Y$ respectively.

\begin{Proposition}\label{preim2}
Let $f:X\to Y$ be a continuous function and $u\in\mathbb{T}_\sigma$.
Then $A\in\mathcal{L}(Y,u)$ implies $A\circ f\in\mathcal{L}(X,u)$.
 \end{Proposition}

{\em Proof.} Let $A\in Q^Y$ be defined by a $u$-family $F$ in $\mathcal{L}(Y)$. Since the
 preimage function $f^{-1}:P(Y)\to P(X)$ is a morphism from $\mathcal{L}(Y)$ to $\mathcal{L}(X)$ by Lemma \ref{measur}, $A\circ f$ is determined by the $u$-family $f^{-1}(F)$ in $\mathcal{L}(X)$ by Lemma \ref{preim}. Therefore, $A\circ f\in\mathcal{L}(X,u)$.
 \qed
 
Next we briefly discuss the relation of $Q$-IFH in $X$ to the Wadge reducibility $\leq_W$ of $Q$-partitions of $X$ (see Introduction). 

\begin{Corollary}\label{wadgeclos}
If $Q$ is antichain (in particular, $Q=\bar{k}$) then any level of the $Q$-IFH in $X$ is closed downwards under Wadge reducibility.
 \end{Corollary}

{\em Proof.} Since $\leq_Q$ is the equality on $Q$,  $A\leq_WB$ iff $A=B\circ f$ for some continuous function $f$ on $X$. Thus, the assertion is a particular case of Proposition \ref{preim2} when $X=Y$.
 \qed

Corollaries \ref{wadgeclos} and \ref{inclus2} show that  Properties (1,5) of the Wadge hierarchy of sets in the Baire space (see the end of Section \ref{wadgedst}) survive under generalisation to the IFH of $k$-partitions in cb$_0$-spaces (the property (1) survives in the weaker form of being bqo).
If $Q$ is not antichain then the closure under Wadge reducibility does not survive in  general, and if $Q$ is not a finite antichain then the levels of $Q$-IFH may be not bqo under inclusion. 

To keep the properties (1,5), one could modify the definition of the $Q$-IFH by taking the closure $\widehat{\mathcal{L}}(X,u)=\{A\in Q^X\mid\exists B(A\leq_WB\in\mathcal{L}(X,u))\}$ of levels under the Wadge reducibility as the new definition. Then we automatically have the closure property (5). It turns out that also the bqo-modification of property (1) holds under this modification. The next proposition is proved in the same way as Theorem \ref{inclus1} and the corresponding lemmas about inclusions of levels of the $Q$-IFS.  

\begin{Proposition}\label{inclus3}
If  $u\trianglelefteq v$ then $\widehat{\mathcal{L}}(X,u)\subseteq\widehat{\mathcal{L}}(X,v)$.
 \end{Proposition}
 
 \begin{Corollary}\label{inclus4}
If $Q$ is a bqo  then  $(\{\widehat{\mathcal{L}}(X,u)\mid u\in\mathbb{T}_\sigma\};\subseteq)$ is also a bqo. 
 \end{Corollary}

One could conclude that taking $\widehat{\mathcal{L}}(X,u)$ instead of $\mathcal{L}(X,u)$ really  improves the definition of $Q$-IFH but it also has the negative effect: the important preservation property from the next subsection holds for classes $\mathcal{L}(X,u)$ but (probably) not for classes $\widehat{\mathcal{L}}(X,u)$. For this reason we prefer to  keep both modifications which are in fact equivalent for the case of $k$-partitions, as we have just discussed.
 
As we know from  Lemma \ref{redbase1}, most of levels of the Borel hierarchy in $X$ have the $\sigma$-reduction property. By Proposition \ref{red-hier}, this implies the following simpler characterisation of many levels of the $Q$-IFH in $X$.

\begin{Proposition}\label{reduc}
For any cb$_0$-space $X$ and any $u\in\mathbb{T}_\sigma$, red-$\mathcal{L}^1(X,u)=\mathcal{L}^1(X,u)$. If $X$ is zero-dimensional then red-$\mathcal{L}(X,u)=\mathcal{L}(X,u)$ for all $u\in\mathbb{T}_\sigma$. 
 \end{Proposition}
 
Let $\mathbf{\Gamma}$ be a family of pointclasses. Recall from \cite{s13} that a total representation (TR)
$\nu:\mathcal{N}\to\mathbf{\Gamma}(X)$ is a {\em $\mathbf{\Gamma}$-TR} if its
{\em universal set} $U_\nu=\{(a,x)\mid x\in\nu(a)\}$ is in
$\mathbf{\Gamma}(\mathcal{N}\times X)$, and $\nu$ is a {\em principal
$\mathbf{\Gamma}$-TR} if it is a $\mathbf{\Gamma}$-TR and any
$\mathbf{\Gamma}$-TR $\mu:\mathcal{N}\to\mathbf{\Gamma}(X)$ is reducible to $\nu$.  Note that if
$\nu:\mathcal{N}\to\mathbf{\Gamma}(X)$ is principal then it is a surjection
and that $\mathbf{\Gamma}(X)$ has at most one principal TR, up to equivalence. According to Theorem 5.2 in \cite{s13}, any level of the classical hierarchies of sets in arbitrary cb$_0$-space has a principal TR. 

The notion of principal TR may be naturally extended to $k$-partitions \cite{s16} and even to $Q$-partitions. Namely, a {\em family of $Q$-partition classes} is a family $\{\mathbf{\Gamma}(X)\}_X$ indexed by cb$_0$-spaces such that $\mathbf{\Gamma}(X)\subseteq Q^X$ for each $X$, and $A\circ f\in\mathbf{\Gamma}(X)$ for every continuous function $f:X\to Y$ and every $A\in\mathbf{\Gamma}(Y)$. A TR $\nu:\mathcal{N}\to\mathbf{\Gamma}(X)$ is a {\em $\mathbf{\Gamma}$-TR} if its
{\em universal $Q$-partition} $(a,x)\mapsto\nu(a)(x)$ is in $\mathbf{\Gamma}(\mathcal{N}\times X)$, and $\nu$ is a {\em principal $\mathbf{\Gamma}$-TR} if it is a $\mathbf{\Gamma}$-TR and any $\mathbf{\Gamma}$-TR $\mu:\mathcal{N}\to\mathbf{\Gamma}(X)$ is reducible to $\nu$.  Note that if $\nu:\mathcal{N}\to\mathbf{\Gamma}(X)$ is principal then it is a surjection and that $\mathbf{\Gamma}(X)$ has at most one principal TR, up to equivalence. 

According to Proposition \ref{preim2}, $\{\mathcal{L}(X,u)\}_X$ is a family of $Q$-partition classes, for each $u\in\mathbb{T}_\sigma$.   But the principal TRs of levels of  $Q$-IFH do not always exist (even for the case of sets).
In particular, for $k$-partitions, $k\geq3$, the  principal TRs of levels of natural hierarchies may not exist. E.g., this is the case already for the difference hierarchies of $3$-partitions over the open sets which consists precisely of the classes $\mathcal{L}(X,T)$, $T\in\mathcal{T}(\bar{3})$. 
A reasonable way to construct a principal TR is to represent all $T$-families of open sets that induce a 3-partition; but this can be done straightforwardly only for reducible bases. Thus, the problem is again related to the $\sigma$-reduction property. This also applies to  iterated labeled trees yielding the following sufficient condition which extends Proposition 4.12 in \cite{s13} and other similar results. The proof consists in ``effectivisation'' of the results above related to the $\sigma$-reduction property. 

\begin{Theorem}\label{reduc1}
Let $X$ be a cb$_0$-space. Then  any level $\mathcal{L}^1(X,u)$ has a  principal total representation. If $X$ is zero-dimensional  then any level $\mathcal{L}(X,u)$ has a principal total representation. 
 \end{Theorem}
 
{\em Proof Sketch.} The proof for both assertions is similar, so we consider only the second one. By Theorem 5.2 in \cite{s13} (see also \cite{brat}), any level $\bfSig^0_{1+\alpha}(X)$ in arbitrary cb$_0$-space has a principal TR $\pi_\alpha$. Moreover, the operations of countable union and binary intersection on $\bfSig^0_{1+\alpha}(X)$ have continuous realizers w.r.t. $\pi_\alpha$. The proof of $\sigma$-reduction property for $\bfSig^0_{1+\alpha}(X)$ also ``effectivizes'', i.e., there is a continuous realizer that computes a reduct of a given (by $\calN$-names) sequence of $\bfSig^0_{1+\alpha}(X)$-sets.

For any given $u\in\mathbb{T}_\sigma$ it suffices to find a TR of the reduced $u$-families in $\mathcal{L}(X)$ that induces the desired principal TR of $\mathcal{L}(X,u)$ by Lemma \ref{tfaml}. 
 If $u$ is a singleton term the TR is obvious. Otherwise, a $u$-family in $\mathcal{L}(X)$ has the form $(\{U_{\tau_0}\},\{F_{\tau_0}\})$ where $\{U_{\tau_0}\}$ is a monotone $u'$-family of  ${\mathcal L}^{sh(u)}_0$-sets with $U_\varepsilon=X$ and, for each $\tau_0\in T(u')$, $F_{\tau_0}$ is  a  $t_0(\tau_0)$-family in $\mathcal{L}^{sh(u)}(\tilde{U}_{\tau_0})$. The TR $\pi_{sh(u)}$ induces a TR of all families $\{U_{\tau_0}\}$. Moreover, by the effective version of Lemma \ref{tfam} we obtain a TR of all monotone such families. By the effective version of Lemma \ref{reduce}, we obtain a TR of all reduced such families.

If the term $t_0(\tau_0)$ is singleton, the procedure of Definition \ref{syntax3} is finished. Otherwise, $F_{\tau_0}=(\{U_{\tau_0\tau_1}\},\{F_{\tau_0\tau_1}\})$ and we can continue the computation above and find a TR of all reduced families $\{U_{\tau_0\tau_1}\}$, for any fixed $\tau_0$. Continuing this process,  we find a desired TR of all $u$-families. This TR induces a TR of $\mathcal{L}(X,u)$ by Lemma \ref{tfaml}. A routine calculation shows that it is a principal $u$-TR.
 \qed 

The  Wadge complete elements in levels of $Q$-IFH do not need to exist. We can prove their existence for the $Q$-IFH in the Baire space (which coincides with the Wadge hierarchy by Theorem \ref{qwadgeBaire} below) but this was  already proved in \cite{km17} by different methods. We give a hint to an elementary proof not using deep facts in \cite{km17}.  

\begin{Corollary}\label{wcomp}
For any $Q$, every  level $\mathcal{L}(\calN,u)$ has  a Wadge complete $Q$-partition.  
 \end{Corollary}
 
{\em Proof.} By Theorem \ref{reduc1}, there is a principal TR  $\nu:\calN\to\mathcal{L}(\calN,u)$, hence its
universal $Q$-partition $U_\nu(a,x)=\nu(a)(x)$ is in $\mathcal{L}(\calN\times\calN,u)$. Since $\calN\times\calN$ is homeomorphic to $\calN$, we can think that $U_\nu\in\mathcal{L}(\calN,u)$. Clearly, any element of $\mathcal{L}(\calN,u)$ is Wadge reducible to $U_\nu$ which is thus Wadge complete in $\mathcal{L}(\calN,u)$.
 \qed

\subsection{Preservation property}\label{preserv}

Here we show that all levels of the $Q$-IFH are preserved by continuous open surjections.

With any function $f:X\to Y$ between cb$_0$-spaces we  associate the  function $A\mapsto f[A]$  from $P(X)$ to $P(Y)$ defined  by
 $f[A]=\{y\in Y\mid A\cap f^{-1}(y)\text{ is non-meager in }f^{-1}(y)\}.$
  Its importance stems from Baire-category properties of cb$_0$-spaces recalled in Section \ref{qpolish}. The function  $A\mapsto f[A]$ (known as the existential category quantifier \cite{ke95,chen}) was used e.g. in \cite{sr07,br,s16}; we changed its notation trying to make it more convenient in our context. 
  
The next two lemmas generalize some results from \cite{sr07,br,s16}. Please distinguish $f[A]$ and the image $f(A)$ of $A$ under $f$.

\begin{Lemma}\label{catim}
\begin{enumerate}
\item The function  $A\mapsto f[A]$ is a morphism from $\mathcal{L}(X)$ to $\mathcal{L}(Y)$, and $f[A]\subseteq f(A)$ for each $A\subseteq X$.
\item If $T$ is a well founded tree and $\{U_\tau\}$ is a $T$-family of $\bfSig^0_{1+\alpha}(X)$-sets then $\{f[U_\tau]\}$ is a $T$-family of $\bfSig^0_{1+\alpha}(Y)$-sets, and $\widetilde{f[U_\tau]}\subseteq f[\tilde{U}_\tau]$ for each $\tau\in T$.
\end{enumerate}
 \end{Lemma}
 
{\em Proof.} (1) Let $y\in f[A]$, then $A\cap f^{-1}(y)$ is non-meager in $f^{-1}(y)$. Then $A\cap f^{-1}(y)$ is non-empty, hence $y\in f(A)$ and $f[A]\subseteq f(A)$. In particular, $f[\emptyset]=\emptyset$. 
To show that $f[X]=Y$ we have to check that, for any $y\in Y$, $f^{-1}(y)$ is non-meager in $f^{-1}(y)$, and this  follows from  quasi-Polishness of $f^{-1}(y)$.
The property that $f[\bigcup_nU_n]=\bigcup_nf[U_n]$ for every countable sequence $\{U_n\}$ in $P(X)$ is well known. 
The (non-trivial) fact that $U\in\bfSig^0_{1+\alpha}(X)$ implies $f[U]\in\bfSig^0_{1+\alpha}(Y)$, follows from Proposition \ref{preserve}, see \cite{sr07,br}.

(2) The first assertion follows from (1), so we check the second one. Let $y\in\widetilde{f[U_\tau]}$, i.e. $y\in f[U_\tau]\setminus\bigcup\{f[U_{\tau'}]\mid \tau\sqsubset\tau'\in T\}$. Then $U_\tau\cap f^{-1}(y)$ is non-meager in $f^{-1}(y)$ and, for each $\tau\sqsubset\tau'\in T$, $U_{\tau'}\cap f^{-1}(y)$ is meager in $f^{-1}(y)$. Then $(\bigcup\{U_{\tau'}\mid \tau\sqsubset\tau'\in T\})\cap f^{-1}(y)$ is meager in $f^{-1}(y)$, hence $\tilde{U}_\tau=U_\tau\setminus\bigcup\{U_{\tau'}\mid \tau\sqsubset\tau'\in T\}$ is non-meager in $f^{-1}(y)$, i.e. $y\in f[\tilde{U}_\tau]$.
 \qed

We associate with any $u$-family $F$ in $\mathcal{L}(X)$ the $u$-family $f[F]$ in  $\mathcal{L}(Y)$ by induction as follows: if $u$ is a singleton term (hence $F=\{X\}$)  then we set $f[F]=\{Y\}$; otherwise, $u'$ is an $F$-term and $F=(\{U_\tau\},\{F_\tau\})$ is a $u'$-family in  $\mathcal{L}^{sh(u)}(X)$; we set $f[F]=(\{f[U_\tau]\},\{f[F_\tau]\})$ which is a $u'$-family in $\mathcal{L}^{sh(u)}(Y)$, hence  a $u$-family in $\mathcal{L}(Y)$.

\begin{Lemma}\label{catim1}
Let $u\in\mathbb{T}_{\sigma}$, $A\in Y\to Q$, and $A\circ f\in\mathcal{L}(X,u)$ be determined by a $u$-family $F$ in $\mathcal{L}(X)$. Then $A$ is determined by the $u$-family $f[F]$ in $\mathcal{L}(X)$. 
 \end{Lemma}
 
 {\em Proof.} If $u$ is a singleton term, the assertion is obvious. Otherwise, $u'$ is an $F$-term and the family $F$ has the form $(\{U_{\tau_0}\},\{U_{\tau_0\tau_1}\},\ldots)$, so $f[F]$ has the form $(\{f[U_{\tau_0}]\},\{f[U_{\tau_0\tau_1}]\},\ldots)$.
We have to show that $A$ is determined by $f[F]$, i.e. for each $y\in Y$,
$A(y)=q(\tau_0,\ldots,\tau_m)$, for every terminating component $\widetilde{f[U_{\tau_0\cdots\tau_m}]}$ of $f(F)$ containing $y$. 
Note that such a component always exists.  

For any given $y\in Y$ and any  such component $\widetilde{f[U_{\tau_0\cdots\tau_m}]}$ we have $y\in f[\tilde{U}_{\tau_0\cdots\tau_m}]$ by Lemma \ref{catim}(2), so $y=f(x)$ for some $x\in\tilde{U}_{\tau_0\cdots\tau_m}$. Thus, $A(y)=(A\circ f)(x)=q(\tau_0,\ldots,\tau_m)$.
 \qed

As an immediate corollary of Lemmas \ref{catim1} and \ref{preim} we obtain the following preservation property for levels of the $Q$-IFH.

\begin{Theorem}\label{qpol}
Let $\mathcal{L}(X),\mathcal{L}(Y)$ be  Borel bases in cb$_0$-spaces  $X,Y$ respectively,  $f:X\to Y$  a continuous open surjection, $A:Y\to Q$, and $u\in\mathbb{T}_{\sigma}$. Then $A\circ f\in\mathcal{L}(X,u)$ iff $A\in\mathcal{L}(Y,u)$.
 \end{Theorem}

{\em Proof.} Let $A\in\mathcal{L}(Y,u)$, then $A$ is determined by a $u$-family $F$ in $\mathcal{L}(Y)$. By Lemma \ref{preim}, $A\circ f\in\mathcal{L}(X,u)$. Conversely, let $A\circ f\in\mathcal{L}(X,u)$, then $A\circ f$ is determined by a $u$-family $F$ in $\mathcal{L}(X)$. By Lemma \ref{catim1}, $A$ is determined by the $u$-family $f[F]$ in $\mathcal{L}(Y)$, hence  $A\in\mathcal{L}(Y,u)$.
 \qed

\subsection{Inheritance of Hausdorff-Kuratowski-type theorems}\label{hktype}

Here we apply the preservation theorem to show that some versions of the  Hausdorff-Kuratowski  theorem (which we call HK-type theorems for short) are inherited by the continuous open images. 

Recall  that the Hausdorff theorem in a space $X$ says that $\bigcup_{\beta<\omega_1}\bfSig^{-1,1}_\beta(X)=\bfDelta^0_{2}(X)$. The difference hierarchy $\{\bfSig^{-1,1}_\beta(X)\}$ over the open sets in $X$ is usually defined using a difference operator on the transfinite sequences of open sets (see e.g. \cite{ke95,s13}). Since in this paper we promote using labeled trees instead of ordinals, we note that  levels $\bfSig^{-1,1}_\beta(X)$ are easily characterised using $\bar{2}$-labeled trees in $\mathcal{T}(\bar{2})$ (see the beginning of Section \ref{trees}). Namely, by Proposition 4.9 in \cite{s13}, there is a tree $T_\beta\in\mathcal{T}(\bar{2})$ such that $\bfSig^{-1,1}_\beta(X)=\mathcal{L}(X,T_\beta)$, and  any $T\in\mathcal{T}(\bar{2})$ is $\trianglelefteq$-equivalent to one of $T_\beta,\bar{T}_\beta$, where $u\mapsto\bar{u}$ is the automorphism induced by $i\mapsto1-i$, see Lemma \ref{auto}. Thus, the Hausdorff theorem for $X$ may be written as $\bigcup\{\mathcal{L}(X,T)\mid T\in\mathcal{T}(\bar{2})\}=\bfDelta^0_{2}(X)$ (in this subsection it is more convenient to work with labeled trees rather that with terms, see Lemma \ref{hqo}).

The Kuratowski theorem extends the Hausdorff theorem to any successor level of the Borel hierarchy in $X$ (see Section \ref{qpolish} for the formulation of this theorem for quasi-Polish spaces). The Kuratowski theorem has a reformulation in terms of  $\bar{2}$-labeled trees in just the same way as for the Hausdorff theorem. Namely, the tree form of the Hausdorff-Kuratowski theorem in $X$  looks like $\bigcup\{\mathcal{L}(X,T)\mid T\in\mathcal{T}_\alpha(\mathcal{T}(\bar{2}))\}=\bfDelta^0_{1+\alpha+1}(X)$ for each $\alpha<\omega_1$, where some notation from the end of Section \ref{trees} is used; in particular, $\mathcal{T}_\alpha\circ\mathcal{T}=\mathcal{T}_{\alpha+1}$.

The tree form of the HK-theorem readily extends to $Q$-partitions which yields our first example of inheritance of the HK-type theorems. We say that a cb$_0$-space $X$ {\em satisfies the HK-theorem for $Q$-partitions in level $1+\alpha+1<\omega_1$}, iff  $\bigcup\{\mathcal{L}(X,T)\mid T\in\mathcal{T}_{\alpha+1}(Q)\}=\bfDelta^0_{1+\alpha+1}(Q^X)$.
We define the qo $\leq_{co}$ on cb$_0$-spaces by: $Y\leq_{co}X$ iff there is a continuous open surjection from $X$ onto $Y$.

\begin{Theorem}\label{hktype2}
If a cb$_0$-space $X$ satisfies the HK-theorem for $Q$-partitions in level $1+\alpha+1<\omega_1$, then so does every space $Y\leq_{co}X$.
\end{Theorem}

{\em Proof.} Since the inclusion $\bigcup\{\mathcal{L}(X,T)\mid T\in\mathcal{T}_{\alpha+1}(Q)\}\subseteq\bfDelta^0_{1+\alpha+1}(Q^X)$ is easy, we check only the opposite inclusion.  Let $A\in\bfDelta^0_{1+\alpha+1}(Q^Y)$ and let $f:X\to Y$ be a continuous open surjection. Then $A\circ f\in\bfDelta^0_{1+\alpha+1}(Q^X)$, hence $A\circ f\in\mathcal{L}(X,T)$ for some $T\in\mathcal{T}_{\alpha+1}(Q)$. By Theorem \ref{qpol}, $A\in\mathcal{L}(Y,T)$. 
 \qed
 
Our second example is concerned with a version of HK-theorem for limit levels of the Borel hierarchy. The problem of finding a construction principle for the $\bfDelta^0_\lambda$-subsets of the Baire
space in the case that $\lambda$ is a positive limit countable ordinal was posed long ago by Luzin and resolved in \cite{wad84} as an important step to the complete description of the Wadge hierarchy. We state the inheritance property for an extension of this result from sets to $Q$-partitions. We say that a cb$_0$-space $X$ {\em satisfies the Wadge property for $Q$-partitions in a limit level $\lambda<\omega_1$}, iff  $\bigcup\{\mathcal{L}(X,T)\mid T\in\mathcal{T}_\lambda(Q)\}=\bfDelta^0_{\lambda}(Q^X)$.

The next result is proved in just the same way as the previous theorem.
 
\begin{Theorem}\label{hktype1}
If a cb$_0$-space $X$ satisfies the Wadge property for $Q$-partitions in a limit level $\lambda<\omega_1$, then so does every space $Y\leq_{co}X$.
\end{Theorem}

\subsection{Characterizing $Q$-Wadge hierarchy in the Baire space}\label{qwadgeb}

Here we show that the $Q$-IFH in the Baire space coincides with the Wadge hierarchy of $Q$-partitions.

The structure of Wadge degrees of Borel measurable $Q$-partitions of $\calN$ was characterised in \cite{km17} (see  Proposition \ref{char1} in Section \ref{trees}). In  particular,  a set-theoretic characterisation of the non-self-dual levels of the $Q$-Wadge hierarchy (with levels $\mathcal{W}(\calN,T)$ from Introduction) was provided (see Lemma 3.16 and its extensions in  \cite{km17}), by defining classes $\Sigma_T$ of $Q$-partitions using set-theoretic operations and showing that $\mathcal{W}(\calN,T)=\widehat{\Sigma}_T$ for each $T\in\mathcal{T}_{\omega_1}(Q)$. 

The definition of $\Sigma_T$ in \cite{km17} uses special features of the Baire space and looks a bit different from our general definition of levels of the $Q$-IFH. The main result of this subsection shows that these classes for the Baire space coincide. For the reader's convenience, we cite necessary notions and results from \cite{km17} (see also \cite{ks19}). 

Any non-empty closed set $C$ in $\calN$ and any $Q$-partition $A:C\to Q$ induce a $Q$-partition $\hat{A}:\calN\to Q$ obtained by composing $A$ with the canonical retraction from $\calN$ onto $C$ (abusing notation, $A$ and $\hat{A}$ are often identified). Similarly, any $A:U\to Q$, where $U$ is a non-empty open set in $\calN$, may be identified with some $\hat{A}:\calN\to Q$ (see Observations 3.5 and 3.6 in \cite{km17}). We recall (in the slightly different from \cite{km17} notation of Section \ref{trees}) the definition of classes $\Sigma_T$ (in fact, we define $\Sigma_u$ for $u\in\mathbb{T}_\sigma$, where $T=T(u)$, see Lemma \ref{hqo}, cf. Definition 3.7 and its extensions in \cite{km17}).

\begin{Definition}\label{sigmau}
\begin{enumerate}
\item $\Sigma_q=\{\lambda x.q\}$.
\item $\Sigma_{s_\alpha(u)}$ consists of $A\circ g$ where $A\in\Sigma_u$ and $g$ is a $\bfSig^0_{1+\omega^\alpha}$-measurable function on $\calN$.
\item $\Sigma_{F_q(u_0,\ldots)}$ consists of $A\in Q^\calN$ such that for some pairwise disjoint non-empty open sets $U_0,U_1,\ldots$ not exhausting $\calN$ we have: $A|_V=\lambda v.q$ (where $V=\calN\setminus\bigcup_iU_i$) and $A|_{U_i}\in\Sigma_{u_i}$ for all $i\geq0$. 
\item $\Sigma_{F_\alpha(u_0,\ldots)}$ consists of $A\in Q^\calN$ such that for some pairwise disjoint non-empty open sets $U_1,U_2,\ldots$ not exhausting $\calN$ we have: $A|_V\in\Sigma_{s_\alpha(u_0)}$ (where $V=\calN\setminus\bigcup_{i\geq1}U_i$) and $A|_{U_i}\in\Sigma_{u_i}$ for all $i\geq1$.
\end{enumerate}
 \end{Definition}

Let $\con\colon\calN\to\calN$ be a function with $\con\circ\con=\con$.
We say that a function $f\colon\calN\to\calN$ is {\em $\con$-conciliatory} if, for any $x,y\in\calN$, $\con(x)=\con(y)$ implies $\con(f(x))=\con(f(y))$.
Similarly, a function $A\colon\calN\to Q$ is {\em $\con$-conciliatory} if, for any $x,y\in\calN$, $\con(x)=\con(y)$ implies $A(x)=A(y)$.
We say that $f,g\colon\calN\to\calN$ are {\em $\con$-equivalent} (written $f\equiv_\con g$) if $\con\circ f=\con\circ g$.

In \cite{km17} the following basic fact was established: For any countable ordinal $\alpha$, there is a $\mathbf{\Sigma}^0_{1+\alpha}$-measurable $\con$-conciliatory function $\mathcal{U}_\alpha\colon\calN\to\calN$ which is {\em universal}; that is, for every $\mathbf{\Sigma}^0_{1+\alpha}$-measurable function $f\colon\calN\to\calN$, there is a continuous function $g\colon\calN\to\calN$ such that $f$ is $\con$-equivalent to $\U_\alpha\circ g$. It was also shown that every $\sigma$-join-irreducible Borel function $A\colon\calN\to Q$ is Wadge equivalent to a $\con$-conciliatory function.
In fact, for any  $u\in\mathbb{T}_{\sigma}$ there is a $\Sigma_u$-complete $\con$-conciliatory function $\mu(u)\colon\calN\to Q$ defined as follows: $\mu(q)=\lambda x.q$; $\mu(s_\alpha(u))=\mu(u)\circ\U_{\omega^\alpha}$; $\mu(F_q(u_0,\ldots))=\mu(q)\to(\mu(u_0)\sqcup\cdots)$; $\mu(F_\alpha(u_0,\ldots))=\mu(s_\alpha(u_0))\to(\mu(u_1)\sqcup\cdots)$.

As usual, let $\mathcal{L}(\calN)$ denote the Borel base in $\calN$. Since this base is reducible, in the proof below we always assume families in $\mathcal{L}(\calN)$ to be reduced (see Proposition \ref{red-hier}). 

\begin{Theorem}\label{qwadgeBaire}
In the Baire space, the $Q$-IFH  coincides with the Wadge hierarchy of $Q$-partitions, i.e. $\Sigma_u=\mathcal{L}(\calN,u)$ for each $u\in\mathbb{T}_{\sigma}$.
 \end{Theorem}
 
{\em Proof.} The equality $\Sigma_q=\mathcal{L}(\calN,q)$ for $q\in Q$ is obvious. To prove $\Sigma_{s_\alpha(u)}=\mathcal{L}(\calN,s_\alpha(u))$,  note that we have $\Sigma_u=\mathcal{L}(\calN,u)$ by induction and that $\mathcal{L}(\calN,s_\alpha(u))=\mathcal{L}^{\omega^\alpha}(\calN,u)$ by Lemma \ref{shift}(1). Let $A\circ g\in\Sigma_{s_\alpha(u)}$ where $A\in\Sigma_u=\mathcal{L}(\calN,u)$ and $g$ is $\mathcal{L}^{\omega^\alpha}(\calN)$-measurable. By Lemmas \ref{measur} and \ref{preim}, $A\circ g\in\mathcal{L}^{\omega^\alpha}(\calN,u)$, as desired. Conversely, let $A\in\mathcal{L}^{\omega^\alpha}(\calN,u)$. By the remarks before the theorem, $\mu(s_\alpha(u))=\mu(u)\circ\U_{\omega^\alpha}$ is Wadge complete in $\mathcal{L}^{\omega^\alpha}(\calN,u)$, hence $A=(\mu(u)\circ\U_{\omega^\alpha})\circ f$ for some continuous function $f$ on $\calN$. Then $A=\mu(u)\circ(\U_{\omega^\alpha}\circ f)$, $\mu(u)\in\mathcal{L}(\calN,u)$, and $\U_{\omega^\alpha}\circ f$ is $\mathcal{L}^{\omega^\alpha}(\calN)$-measurable. Thus, $A\in\Sigma_{s_\alpha(u)}$.

In proving the equality $\Sigma_{F_q(u_0,\ldots)}=\mathcal{L}(\calN,F_q(u_0,\ldots))$, by induction we can assume that $\Sigma_{u_i}=\mathcal{L}(\calN,u_i)$ for each $i\geq0$. Let $A\in\Sigma_{F_q(u_0,\ldots)}$, then for some pairwise disjoint non-empty open sets $U_0,U_1,\ldots$ not exhausting $\calN$ we have: $A|_V=\lambda v.q$ (where $V=\calN\setminus\bigcup_iU_i$) and $A|_{U_i}\in\Sigma_{u_i}$ for all $i\geq0$. By induction,  $A|_{U_i}\in\mathcal{L}(\calN,u_i)$ for all $i\geq0$. By Lemma \ref{shift}(4), $A\in\mathcal{L}(\calN,F_q(u_0,\ldots))$. The converse inclusion follows from Lemma \ref{shift}(3) and Definition \ref{sigmau}(3). The case of $F_\alpha$-term is considered similarly.
 \qed

\subsection{Infinitary fine hierarchies in quasi-Polish spaces}\label{qfineq}

Here we summarise some properties of the $Q$-IFH in quasi-Polish spaces. For any quasi-Polish space $X$ we fix a continuous open surjection $\xi$ from $\calN$ onto $X$ (Proposition \ref{adm}). 
First we give the characterisation of the Wadge hierarchy of $k$-partitions announced in Introduction (for $k=2$ this of course yields a characterisation of the Wadge hierarchy of sets).

\begin{Theorem}\label{qwadgeqp}
Let $X$ be a quasi-Polish space, $Q=\bar{k}$, and $T\in\mathcal{T}_{\omega_1}(Q)$. Then $\mathcal{W}(X,T)=\mathcal{L}(X,T)$.
 \end{Theorem}
 
 {\em Proof.} Let $\xi:\calN\to X$ be a continuous open surjection. By  Theorem \ref{qwadgeBaire} and Proposition \ref{char1}, $\mathcal{W}(\calN,T)=\Sigma_T=\mathcal{L}(\calN,T)$. By Theorem \ref{qpol}, for any $A:X\to Q$ we have: $A\in\mathcal{W}(X,T)$ iff $A\circ\xi\in\mathcal{L}(\calN,T)$ iff $A\in\mathcal{L}(X,T)$.
 \qed

Next we show that the HK-type theorems hold in any quasi-Polish space, which extends some known results.
From Proposition \ref{adm} we know that $X$ is a quasi-Polish space iff $X\leq_{co}\mathcal{N}$. This together with Theorems \ref{hktype2} and \ref{hktype1} implies the following.

\begin{Theorem}\label{hausp}
Every quasi-Polish space satisfies the HK-theorem for $Q$-partitions in any successor level $1+\alpha+1<\omega_1$ of the $Q$-IFH, and also the Wadge property for $Q$-partitions in  any limit level $\lambda<\omega_1$ of the $Q$-IFH.
 \end{Theorem}  

Next we make some remarks on which properties of the Wadge hierarchy in the Baire space (see end of Section \ref{wadgedst}) hold in arbitrary quasi-Polish spaces. Property (1) holds for the hierarchies of sets and of $k$-partitions (for $k\geq3$ in the weakened bqo-form); the non-collapse property (2) does not automatically hold and requires an additional investigation in any concrete space; property (3) fails in most of natural spaces; property (4) holds in arbitrary quasi-Polish space (note that this property is in fact an HK-type theorem); property (5) holds for the hierarchies of sets and of $k$-partitions; property (6) does not automatically hold and requires an additional investigation in any concrete space. % property (7) holds for many $\bfPi$-levels (as we show below).

\section{Effective Wadge hierarchy}\label{effect}

Here we briefly discuss  effective versions of some notions and results described so far. For a detailed presentation of the effective versions see \cite{s19}.

For Theoretical  Computer Science and Computable Analysis  an effective  DST for reasonable classes of effective spaces is especially relevant. A lot of useful work in this direction was done in Computability Theory but mostly for the discrete space $\omega$, the Baire space $\mathcal{N}$, and some of their  relatives \cite{ro67,mo09}. 
Effective versions of the classical Borel, Hausdorff and Luzin hierarchies are naturally defined for every effective space (see e.g.\cite{s06,s15}) but, as also in the classical case, they behave well only for spaces of special kind. 

By {\em  effectivization of a cb$_0$-space $X$} we mean a numbering $\beta:\omega\to P(X)$ of a base in $X$ such that there is a uniform sequence $\{A_{ij}\}$ of c.e. sets with
$\beta_i\cap\beta_ j=\bigcup\beta(A_{ij})$, where $\beta(A)$ is the image of $A$ under $\beta$. The numbering $\beta$ is called an {\em effective base of $X$} while the pair $(X,\beta)$ is called an {\em effective space}. We simplify $(X,\beta)$ to $X$ if $\beta$ is clear from the context. The \emph{effectively open sets} in $X$ are the sets~$\bigcup_{i\in W}\beta(i)$, for some c.e.~set~$W\subseteq\omega$. The standard numbering $\{W_n\}$ of c.e. sets \cite{ro67} induces a numbering of the effectively open sets. The notion of effective space allows to define e.g. computable and effectively open functions between such spaces \cite{wei00}.

Recently, a convincing version of a computable quasi-Polish space (CQP-space for short) was suggested in \cite{br1,hs}. Effective versions of some classical facts (e.g., of the Hausdorff theorem) were established in \cite{s15} for CQP-spaces. By a {\em computable quasi-Polish space}  we mean an effective space  $(X,\beta)$ such that there exists a computable effective open surjection $\xi:\mathcal{N}\to X$ from the Baire space onto $(X,\beta)$. As shown in \cite{s15,br1,hs},  CQP-spaces do satisfy effective versions of several important properties of quasi-Polish spaces. E.g. they  subsume  computable Polish spaces and  computable domains and satisfy the effective Hausdorff and Suslin theorems. The class of CQP-spaces includes most  of cb$_0$-spaces considered in the literature.

Let  $\{\Sigma^0_{1+n}(X)\}_{n<\omega}$ be the effective Borel hierarchy and  $\{D_n(\Sigma^0_m(X))\}_n$ be the effective Hausdorff difference hierarchy over $\Sigma^0_m(X)$  in arbitrary effective space $X$. Another popular notation for levels of the difference hierarchy is $\Sigma^{-1,m}_n=D_n(\Sigma^0_m(X))$, with $\Sigma^{-1,1}$ usually simplified to $\Sigma^{-1}$.  Let also $\{\Sigma^1_{1+n}(X)\}$ be the effective Luzin hierarchy. We do not repeat the standard definitions (which may be found e.g. in \cite{s06,s15})  but mention that the definitions yield also standard numberings of all levels of the hierarchies, so we can speak e.g. about uniform sequences of sets in a given level. E.g., $\Sigma^0_1(X)$ is the class of effectively open sets in $X$,   $\Sigma^{-1}_2(X)$ is the class of differences of $\Sigma^0_1(X)$-sets, and $\Sigma^0_2(X)$ is the class of effective countable unions of $\Sigma^{-1}_2(X)$-sets.  

Levels
of the effective hierarchies are denoted in the same manner as levels of the
corresponding classical hierarchies, using the lightface letters
$\Sigma,\Pi,\Delta$ instead of the boldface
$\bf{\Sigma},\bf{\Pi},\bf{\Delta}$ used for the classical
hierarchies. The boldface classes may be considered as ``limits''  of the corresponding light-face levels (where the ``limit'' is obtained by taking the union of the corresponding relativised light-face levels, for all oracles). Thus, the effective hierarchies not only refine but also generalise the classical ones.

In \cite{s19}  we developed an effective Wadge hierarchy (including the hierarchy of $k$-partitions) in effective spaces which subsumes the effective Borel and  Hausdorff hierarchies (as well as many others) and is in a sense the finest possible hierarchy of effective Borel sets. By effective Wadge hierarchy in a given effective space we mean the fine hierarchy over $\{\Sigma^0_{1+n}(X)\}_{n<\omega}$ (see e.g. \cite{s08} for a survey). Roughly speaking, the FH is a finitary version of the IFH where one  uses $\omega$ instead of $\omega_1$ and finite trees instead of well founded trees. The finitary analogue of $\mathcal{T}_{\omega_1}(Q)$ is denoted as  $\mathcal{T}_{\omega}(Q)$ and considered in \cite{s19} only for $Q=\bar{k}$.

E.g.,  a {\em base in a set $X$} is now a sequence $\mathcal{L}=\{\mathcal{L}_n\}_{n<\omega}$ of subsets of $P(X)$ such that any $\mathcal{L}_n$ is closed under union and intersection, contains $\emptyset,X$ and satisfies $\mathcal{L}_n\cup\check{\mathcal{L}}_n\subseteq\mathcal{L}_{n+1}$. The {\em effective Borel bases} $\mathcal{L}(X)=\{\Sigma^0_{1+n}(X)\}$ in  effective spaces $X$ are especially relevant.
The (finitary) FH of sets over the base $\mathcal{L}$ is now a sequence $\{\mathcal{S}_\alpha\}_{\alpha<\varepsilon_0}$, $\varepsilon_0=sup\{\omega,\omega^\omega,\omega^{\omega^\omega},\dots\}$, of subsets of $P(X)$ constructed from the sets in levels of the base in $X$ by induction on $\alpha$ using suitable set operations. 

The FH over the effective Borel base in $X$ will be denoted by  $\{\Sigma_\alpha(X)\}_{\alpha<\varepsilon_0}$ and called the {\em effective Wadge hierarchy in $X$}. 
Denote the corresponding boldface sequence by $\{\mathbf{S}_\alpha(\mathcal{N})\}_{\alpha<\varepsilon_0}$. The sequence $\{\mathbf{S}_\alpha(\mathcal{N})\}_{\alpha<\varepsilon_0}$ forms a small but important fragment of the classical Wadge hierarchy in the Baire space studied e.g. in \cite{s03}. In the classical Wadge hierarchy these pointclasses have of course different notations.  It is not hard to show that ${\mathbf S}_\alpha={\mathbf
\Sigma}_{f(\alpha)}$ for each $\alpha<\varepsilon_0$ where
$f:\varepsilon_0\rightarrow\upsilon$ is the monotone function defined by
induction as follows: \( f(0)=0 \) and
 $$f(\omega^{\alpha_1}\cdot k_{1}+\omega^{\alpha_2}\cdot
k_{2}+\ldots )= \omega_{1}^{f(\alpha_1)}\cdot k_{1}+\omega
_{1}^{f(\alpha_2)}\cdot k_{2}+\ldots,$$ for any non-empty sequence
\( \alpha_1>\alpha_2>\ldots  \) of ordinals $<\varepsilon_0$, and
for all \( k_{i}<\omega \) (recall that any positive ordinal
$\alpha<\varepsilon_0$ is uniquely representable in the form
$\alpha=\omega^{\alpha_1}\cdot k_{1}+\omega^{\alpha_2}\cdot
k_{2}+\ldots $).

The finitary FH of $k$-partitions  over the effective Borel base is denoted  as $\{\Sigma(X,T)\}_{T\in\mathcal{T}_{\omega}(\bar{k})}$.
In particular, we show in \cite{s19} that levels of such hierarchies are preserved by the computable effectively open surjections, that if the effective Hausdorff-Kuratowski theorem holds in the Baire space then it holds in every CQP-space,  and we extend the effective Hausdorff theorem for CQP-spaces \cite{s15} to $k$-partitions. We hope that these results (together with those already known) show that the effective DST has already reached the state of maturity.

\section{Future work}\label{con}

Many interesting questions related to this paper remain open even for the case of $k$-partitions $Q=\bar{k}=\{\mathbf{0},\ldots,\mathbf{k-1}\}$. We shorten the signature $\sigma(\bar{k},\omega_1)$ to $\sigma(k)$  (the boldface symbols are used to distinguish the elements of $\bar{k}$ from ordinals $0,\ldots,k-1$). By Proposition \ref{char1}, the quotient-poset of $(\mathbb{T}_{\sigma(k)};\trianglelefteq)$ contains essential information about the Wadge hierarchy of ($\sigma$-join-irreducible) $k$-partitions of the Baire space. But if for $k=2$ most questions about the structure $(\mathbb{T}_{\sigma(k)};\trianglelefteq)$ follow from the results in \cite{wad84}, for $k\geq3$ there is still a lot to do. Below we assume  that $k\geq3$. 

In \cite{ksz09} it was shown that the automorphism group of the quotient-poset of a natural initial segment of $(\mathbb{T}_{\sigma(k)};\trianglelefteq)$ is isomorphic to the symmetric group $\mathbf{S}_k$ on $k$ elements. We guess that this result extends to  the quotient-poset of $(\mathbb{T}_{\sigma(k)};\trianglelefteq)$.

In \cite{ks09} it was shown that the first-order theory of the quotient-poset of a small initial segment of $(\mathbb{T}_{\sigma(k)};\trianglelefteq)$ is computably isomorphic to the first-order arithmetic; this implies that the first-order arithmetic is $m$-reducible to the   first-order theory of the quotient-poset of  $(\mathbb{T}_{\sigma(k)};\trianglelefteq)$. This is in contrast with the quotient-poset of $(\mathbb{T}_{\sigma(2)};\trianglelefteq)$ whose first-order theory is decidable. Also, in \cite{ks09} a complete characterisation of first-order definable relations in the mentioned small initial segment was achieved. This motivates the study of definable relations on the quotient-poset of  $(\mathbb{T}_{\sigma(k)};\trianglelefteq)$; along with first-order definability, the $L_{\omega_1\omega}$-definabilty in this quotient-poset seems especially interesting.

In this paper we hopefully found a convincing set-theoretic definition of  $Q$-Wadge hierarchy in quasi-Polish spaces, restricting our attention to Borel $Q$-partitions. For this the axioms of ZFC suffice.
A major open question  is to extend the results of this paper to a reasonable class beyond the Borel $Q$-partitions (perhaps even to all $Q$-partitions). 
The Wadge hierarchy for arbitrary subsets of the Baire space is well known \cite{vw76} and requires suitable set-theoretic axioms alternative to  ZFC. The definitions of this paper extend straightforwardly (by taking arbitrarily large ordinal $\gamma$ in the signature $\sigma(Q,\gamma)$ in Section \ref{trees}) but beyond the Borel $Q$-partitions proofs could turn out different from those used in this paper. It is currently not clear which set-theoretic axioms should be used.


\begin{thebibliography}{aasa67}

%\bibitem{aj}   Abramsky S.,   Jung, A.: Domain theory. In: {\em Handbook of Logic in Computer Science}, v. 3, Oxford,1994, 1--168.

%\bibitem{ad65}  Addison J.W. The method of alternating chains. In: {\em The theory of models}, Amsterdam, North Holland, 1965, p.1--16.

\bibitem{brat} Brattka V. Effective borel measurability and reducibility of functions. Mathematical Logic
369 Quarterly, 51(1), p. 19--44, 2005. doi:10.1002/malq.200310125.

\bibitem{bh02} Brattka V. and Hertling P. Topological
properties of real number representations. {\em Theoretical
Computer Science}, 284 (2002), 241--257.

\bibitem{bg15}  Becher V. and Grigorieff S. Wadge hardness in Scott spaces and its effectivization, Math. Structures Comput. Sci. 25 (2015), no. 7, 1520--1545.

\bibitem{chen} Chen R. Notes on quasi-Polish spaces, Submitted. arXiv 1809.07440.

\bibitem{ch}   Callard A.,   Hoyrup M. Descriptive complexity of non-Polish spaces. Submitted to STACS-2020.  

%\bibitem{dp94} Davey, B.A., Priestley, H.A.: {\em Introduction to Lattices and Order}. Cambridge, 1994.

\bibitem {br} %[Br13]
de Brecht M.
Quasi-Polish spaces,
\emph{Annals of pure and applied logic}, \textbf{164}, (2013), 356--381.

\bibitem{br1}
 de Brecht M.,   Pauly A.,  Schr\"oder  M. Overt choice.  arXiv 1902.05926v1 [math. LO], 2019.

\bibitem{du01} Duparc J. Wadge hierarchy and Veblen hierarchy, part I. Journal of Symbolic Logic, 66, No 1 (2001), 56--86.

\bibitem{du19} Duparc J., Vuilleumier L. The Wadge order on the Scott domain is not a well-quasi-order. To appear in J. of Symbolic Logic, DOI: https://doi.org/10.1017/jsl.2019.51, Arxiv 1902.09419.  

\bibitem{ems87}  van Engelen F., Miller A., Steel J. Rigid Borel
sets and better quasiorder theory.  Contemporary mathematics,
65 (1987), 199--222.

\bibitem{en89} %[En89]
Engelking R. 
General Topology,
\emph{Heldermann, Berlin} (1989).

\bibitem{he93}
 Hertling P. Topologische Komplexit\"atsgrade von
Funktionen mit endlichem Bild.  Informatik-Berichte 152, 
Fernuniversit\"at Hagen, 1993.

\bibitem{he96}
Hertling P. {\em Unstetigkeitsgrade von Funktionen in der
effektiven   Analysis}. PhD thesis, Fachbereich Informatik,
FernUniversit\"{a}t Hagen, 1996.

\bibitem{hs}   Hoyrup M.,   Rojas C.,  Selivanov  V.,  Stull D. Computability on quasi-Polish spaces. Proc. of DCFS-2019,  LNCS volume 11612, Berlin, Springer, 2019,   P. 171--183.

\bibitem{ik10}  Ikegami D. Games in Set Theory and Logic, PhD Thesis, University of Amsterdam, 2010.

\bibitem{sc10}  Ikegami D.,  Schlicht P.,  Tanaka H. Continuous reducibility for the real line, preprint, submitted, 2012.

\bibitem{ke95} %[Ke95]
Kechris A.S.
Classical Descriptive Set Theory,
\emph{Springer, New York}, (1995).

\bibitem{kls12}  Kechris A.S.,  L\"owe B, and Steel J.R. (eds.)
Wadge degrees and projective ordinals. The Cabal Seminar. Volume
II, Lecture Notes in Logic, vol. 37, Association for Symbolic Logic,
La Jolla, CA; Cambridge University Press, Cambridge, 2012.

\bibitem{km67}  Kuratowski K., Mostowski A.
{\em Set Theory}. North Holland, 1967.

\bibitem{km17}  Kihara T.,  Montalban A. On the structure of the Wadge degrees of bqo-valued Borel functions. Transactions of the American Mathematical Society, 371, No. 11 (2019), 7885--7923.

\bibitem{ks19}  Kihara T.  Selivanov V. Wadge-like degrees of Borel bqo-valued functions. Submitted, Arxiv 1909.10835.

\bibitem{ks09}  Kudinov O.V., Selivanov V.L. A Gandy theorem for abstract structures and applications to first-order definability. Proc. CiE-2009, LNCS 5635, Springer, Berlin, 290--299, 2009. 

\bibitem{ksz09} Kudinov O.V., Selivanov V.L., Zhukov A.V. Definability in the h-quasiorder of labeled forests. Annals of Pure and Applied Logic, 159(3), 318--332, 2009.

%\bibitem{la71}  Laver, R.: On Fra\"isse's order type conjecture. Ann. of Math. 93(2), 89--111 (1971).

%\bibitem{la78}  Laver, R.: Better-quasi-orderings and a class of trees, in Studies in Foundations and Combinatorics,
%Adv. in Math. Supplementary Series 1, Academic Press, New York, 1978, pp. 31--48.

\bibitem{lo12}  Louveau A. Some results in the Wadge hierarchy of Borel sets,
In \cite{kls12}, 2012, pp 47--73.

\bibitem{mo09}  Moschovakis Y.N. {\em Descriptive  Set  Theory},
North Holland, Amsterdam, 2009.

\bibitem{mss15} Motto Ros L., Schlicht P., Selivanov V.L. Wadge-like reducibilities on arbitrary quasi-Polish spaces. Mathematical Structures in Computer Science, 25, Special Issue 08 (2015), 1705--1754.

%\bibitem{nw65}  Nash-Williams, C. St. J. A.: On well-quasi-ordering infinite trees, Proc. Cambridge Philos. Soc., v. 61 (1965), 697--720.

\bibitem{pe15}  Pequignot Y. A Wadge hierarchy for second countable spaces.
Archive for Mathematical Logic 54.5-6 (2015), pp. 659--683. doi:
10.1007/s00153-015-0434-y.

\bibitem{ro67}   Rogers H., jr. {\em  Theory of Recursive Functions and
Effective Computability}. McGraw-Hill, New York, 1967.

%\bibitem{s83}  Selivanov V.L. Hierarchies  of hyperarithmetical  sets and functions. {\em Algebra i Logika,} 22, No 6 (1983), p.666--692 (English translation: {\em Algebra and Logic,} 22 (1983), p.473--491).

%\bibitem{s95} 	 Selivanov V.L. Fine hierarchies and Boolean terms. The Journal of Symbolic Logic, 60, No 1 (1995), 289--317.

\bibitem{s03}
 Selivanov V.L. Wadge degrees of $\omega$-languages of deterministic Turing machines.
{\em Theoretical Informatics and Applications}, 37 (2003), 67--83.

\bibitem{s05} Selivanov V.L. Variations on Wadge reducibility. Siberian Advances in Mathematics, 15, N 3 (2005), 44--80. 

\bibitem{s06}  Selivanov V.L. Towards a descriptive
set theory for domain-like structures. {\em
Theoretical Computer Science}, 365 (2006), 258--282.

\bibitem{s07}  Selivanov V.L. The quotient algebra of labeled
forests modulo h-equivalence.
{\em Algebra and Logic}, 46, N 2 (2007), 120--133.

\bibitem{s07a}  Selivanov V.L. Hierarchies of
${\mathbf\Delta}^0_2$-measurable $k$-partitions.  Math. Logic
Quarterly, 53 (2007), 446--461.

\bibitem{s08} Selivanov V.L. Fine hierarchies and $m$-reducibilities in
theoretical computer science. {\em  Theoretical Computer Science},
405 (2008), 116--163.

\bibitem{s12} Selivanov V.L.
Fine hierarchies via Priestley duality. Annals of Pure and Applied
Logic, 163 (2012) 1075--1107.

\bibitem{s13}
Selivanov V.L.
Total representations.
\emph{Logical Methods in Computer Science}, 9(2) (2013),
p. 1--30.

\bibitem{s15} Selivanov V.L. Towards the effective descriptive set theory. Proc. CiE 2015, LNCS volume 9136, Berlin, Springer, 2015, P. 324--333.

\bibitem{s16} Selivanov V.L. Towards a descriptive theory of cb0-spaces. Mathematical Structures in Computer Science. v. 28 (2017), issue 8, 1553--1580. DOI: http://dx.doi.org/10.1017/S0960129516000177. Earlier version in: ArXiv: 1406.3942v1 [Math.GN] 16 June 2014.

\bibitem{s17} Selivanov V.L. Extending Wadge theory to k-partitions. J. Kari, F. Manea and Ion Petre (eds.) CiE 2017, LNCS 10307, 387--399, Berlin, Springer.

\bibitem{s18} Selivanov V.L. Q-Wadge degrees as free structures, Submitted.

\bibitem{s19} Selivanov V.L. Effective Wadge hierarchy in computable quasi-Polish spaces, submitted. Arxiv 1910.13220.

%\bibitem{sc11} Schlicht, P.: Continuous reducibility for Polish spaces, 2012, submitted.

%\bibitem{si87} Logic and Combinatorics, edited by S.G. Simpson, Contemporary Mathematics 65, American Mathematical Society, Providence, RI, 1987.

\bibitem{si85}  Simpson S.G. Bqo-theory and Fra\"iss\'e conjecture. Chapter 9 of R. Mansfield, G. Weitkamp. Recursive aspects of descriptive set theory. Oxford
University Press, New York, 1985.

\bibitem{sr07}
Saint Raymond J. Preservation of the Borel class under countable-compact-covering mappings. Topology and its Applications 154 (2007), 1714--1725.

\bibitem{ste80}  Steel J. Determinateness and the separation property.
{\em J. Symbol. Logic,} 45 (1980), p.143--146.

\bibitem{vw76}   Van Wesep R. Wadge  degrees  and  descriptive set theory.
 Lecture Notes in Mathematics, 689 (1976), p. 151--170.

\bibitem{wad72}  Wadge W. Degrees of complexity of subsets of the Baire space. {\em Notices AMS}, 1972, A-714.

\bibitem{wad84}  Wadge W.  Reducibility and Determinateness in the
Baire Space. PhD thesis, University of California, Berkely, 1983.

\bibitem{wei00}  Weihrauch K.
Computable Analysis. Springer Verlag, Berlin (2000).





\end{thebibliography}
\end{document}